\documentclass[a4paper,UKenglish]{lipics_md}
\usepackage{microtype}
\usepackage{mathtools}
\usepackage[normalem]{ulem}

\title{Random Inscribed Polytopes Have Similar Radius Functions
       as Poisson--Delaunay Mosaics\footnote{This work is partially
         supported by the Austrian Science Fund (FWF),
         grant no.\ I02979-N35.}}
\titlerunning{Random Inscribed Polytopes}

\author[1]{Herbert Edelsbrunner}
\author[1]{Anton Nikitenko}
\affil[1]{IST Austria (Institute of Science and Technology Austria),
  Am Campus 1, \\ 3400 Klosterneuburg, Austria,
  \texttt{edels@ist.ac.at}, \texttt{anton.nikitenko@ist.ac.at}}
\authorrunning{H. Edelsbrunner and A. Nikitenko}

\mathsubjclass{60D05 Geometric probability and stochastic geometry,
  68U05 Computer graphics; computational geometry.} 

\newcommand{\mm}[1] {\ifmmode{#1}\else{\mbox{\(#1\)}}\fi}

\newcommand{\ignore}[1]{}

\newcommand{\ourproof}{\begin{proof}}
\newcommand{\eop}{\end{proof}}  %

\newcommand{\Aspace}        {\mm{{\mathbb A}}}
\newcommand{\Bspace}        {\mm{{\mathbb B}}}
\newcommand{\Rspace}        {\mm{{\mathbb R}}}
\newcommand{\Sspace}        {\mm{{\mathbb S}}}
\newcommand{\EE}            {\mm{{\mathbb E}}}
\newcommand{\Expected}[1]   {\mm{{\EE}{[{#1}]}}}
\newcommand{\PP}            {\mm{{\mathbb P}}}
\newcommand{\Probable}[1]   {\mm{{\PP}{[{#1}]}}}
\newcommand{\PemptyOnly}    {\mm{{\PP}_{\emptyset}}}
\newcommand{\Pempty}[1]     {\mm{{\PP}_{\emptyset}{({#1})}}}
\newcommand{\One}[1]        {\mm{{{\bf 1}_{#1}}}}

\newcommand{\Gama}[1]       {\mm{{\Gamma}{({#1})}}}
\newcommand{\iGama}[2]      {\mm{{\gamma}{({#2}; {#1})}}}
\newcommand{\Beta}[2]       {\mm{{B}{({#1},{#2})}}}
\newcommand{\iBeta}[3]      {\mm{{B}_{#1}{({#2},{#3})}}}
\newcommand{\LGrass}[2]     {\mm{{\cal L}_{#1}^{#2}}}
\newcommand{\Delaunay}[2]   {\mm{{\rm Del}_{#1}{{#2}}}}
\newcommand{\Voronoi}[1]    {\mm{{\rm Vor}{({#1})}}}
\newcommand{\Capp}[2]       {\mm{{\rm Cap}_{#1}{({#2})}}}
\newcommand{\capp}[1]       {\mm{{\rm cap}{({#1})}}}
\newcommand{\ecapp}[1]      {\mm{{\rm cap}_{\emptyset}{({#1})}}}
\newcommand{\plane}         {\mm{P}}
\newcommand{\Rfun}          {\mm{{\cal R}}}
\newcommand{\density}       {\mm{\rho}}
\newcommand{\NRad}          {\mm{\bar{\eta}}}
\newcommand{\GRad}          {\mm{\eta}}
\newcommand{\GRadPr}        {\mm{\zeta}}
\newcommand{\ERad}          {\mm{r}}
\newcommand{\littleoh}[1]   {\mm{o{\left({#1}\right)}}}
\newcommand{\littleomega}[1]{\mm{\omega{\left({#1}\right)}}}
\newcommand{\bigoh}[1]      {\mm{O{({#1})}}}
\newcommand{\bigtheta}[1]   {\mm{\Theta{({#1})}}}

\newcommand{\uuu}           {\mm{{\bf u}}}
\newcommand{\xxx}           {\mm{{\bf x}}}
\newcommand{\yyy}           {\mm{{\bf y}}}

\newcommand{\Ccon}[3]       {\mm{{C}_{{#1},{#2}}^{#3}}}
\newcommand{\Dcon}[2]       {\mm{{D}_{#1}^{#2}}}
\newcommand{\ccon}[3]       {\mm{{c}_{{#1},{#2}}^{#3}}}
\newcommand{\dcon}[2]       {\mm{{d}_{#1}^{#2}}}

\newcommand{\EAux}[3]       {\mm{{E}_{{#1},{#2}}^{#3}}}
\newcommand{\Fraction}[1]   {\mm{F}{({#1})}}
\newcommand{\Volume}[1]     {\mm{\rm Vol}{({#1})}}
\newcommand{\Area}[1]       {\mm{\rm Area}{({#1})}}
\newcommand{\Length}[1]     {\mm{\rm Length}{({#1})}}

\newcommand{\card}[1]       {\mm{|{#1}|}}
\newcommand{\dime}[1]       {\mm{\rm dim\,}{#1}}
\newcommand{\conv}[1]       {\mm{\rm conv\,}{#1}}
\newcommand{\interior}[1]   {\mm{\rm int\,}{#1}}
\newcommand{\boundary}[1]   {\mm{\rm bd\,}{#1}}
\newcommand{\diff}          {\mm{\rm \,d}}
\newcommand{\norm}[1]       {\mm{\|{#1}\|}}
\newcommand{\Edist}[2]      {\mm{\|{#1}-{#2}\|}}
\newcommand{\GG}            {\mm{{d}}}
\newcommand{\Gdist}[2]      {\mm{{\GG}{({#1},{#2})}}}
\newcommand{\dd}            {\mm{\delta}}
\newcommand{\ee}            {\mm{\varepsilon}}

\newcommand{\const}         {\mm{{\mathrm{const}}}}

\newcommand{\ourparagraph}[1]  {\vspace{0.1in} \noindent \textbf{#1}}
\newcommand{\Skip}[1]       {}

\begin{document}
\maketitle

\begin{abstract}
  Using the geodesic distance on the $n$-dimensional sphere, we study the
  expected radius function of the Delaunay mosaic of a random set of points.
  Specifically, we consider the partition of the mosaic into
  intervals of the radius function and determine the expected number
  of intervals whose radii are less than or equal to a given threshold.
  Assuming the points are not contained in a hemisphere, the Delaunay mosaic
  is isomorphic to the boundary complex of the convex hull in $\Rspace^{n+1}$,
  so we also get the expected number of faces of a random inscribed polytope.
  We find that the expectations are essentially the same as for the
  Poisson--Delaunay mosaic in $n$-dimensional Euclidean space.
  As proved by Antonelli and collaborators \cite{Antonelli}, an
  orthant section of the $n$-sphere is isometric to the standard
  $n$-simplex equipped with the Fisher information metric.
  It follows that the latter space has similar stochastic properties
  as the $n$-dimensional Euclidean space.
  Our results are therefore relevant in information geometry
  and in population genetics.
\end{abstract}

%\newpage
%%%%%%%%%%%%%%%%%%%%%%%%%%%%%%%%%%%%%%%%%%%%%%%%%%%%%%%%%%%%%%%%%%%%%%%%%%
%%%%%%%%%%%%%%%%%%%%%%%%%%%%%%%%%%%%%%%%%%%%%%%%%%%%%%%%%%%%%%%%%%%%%%%%%%
\section{Introduction}
\label{sec:1}
%%%%%%%%%%%%%%%%%%%%%%%%%%%%%%%%%%%%%%%%%%%%%%%%%%%%%%%%%%%%%%%%%%%%%%%%%%
%%%%%%%%%%%%%%%%%%%%%%%%%%%%%%%%%%%%%%%%%%%%%%%%%%%%%%%%%%%%%%%%%%%%%%%%%%

Letting $X$ be a Poisson point process
in $\Rspace^n$, the expected sizes of the Voronoi tessellation and,
equivalently, of the dual Delaunay mosaic are reasonably well understood.
The starting point for this paper is the question how these expectations
change when we pick the points on the $n$-dimensional sphere, $\Sspace^n$.
Perhaps surprisingly, the difference is very small.
Even the partitions of the Delaunay mosaics into the intervals
of the respective radius functions are barely distinguishable.

%%%%%%%%%%%%%%%%%%%%%%%%%%%%%%%%%%%%%%%%%%%%%%%%%%%%%%%%%%%%%%%%%%%%%%%%%%
\ourparagraph{Motivation.}
%%%%%%%%%%%%%%%%%%%%%%%%%%%%%%%%%%%%%%%%%%%%%%%%%%%%%%%%%%%%%%%%%%%%%%%%%%
Our reason for comparing random sets in the Euclidean space and on the sphere
is the Fisher information metric, which measures the dissimilarity between
discrete probability distributions.
Write $\xxx = (x_0, x_1, \ldots, x_n)$ and $\yyy = (y_0, y_1, \ldots, y_n)$
for two such distributions,
with $\sum_{i=0}^n x_i = \sum_{i=0}^n y_i = 1$ and $x_i, y_i \geq 0$
for all $i$,
and note that $\xxx$ and $\yyy$ are points of the $n$-dimensional standard simplex,
$\Delta^n$.
Letting $\gamma \colon [0,1] \to \Delta^n$ be a smooth curve connecting
$\xxx = \gamma(0)$ to $\yyy = \gamma(1)$, we define its \emph{length} as
\begin{align}
  \Length{\gamma}  &=  \int\displaylimits_{t=0}^1 \sqrt{ \tfrac{1}{2}
    \sum\nolimits_{i=0}^n \frac{\dot\gamma_i(t)^2}{\gamma_i(t)}} \diff t ,
\end{align}
in which $\gamma_i (t)$ and $\dot \gamma_i (t)$ are the $i$-th components
of the curve and its velocity vector.
The \emph{Fisher information metric} assigns the length
of the shortest connecting path
to the pair $\xxx, \yyy$;
see \cite[Section 2.2]{AmNa00} as well as \cite[Section I.4]{Aki79},
where this metric is referred to as the Shahshahani metric.
This way of measuring distance is fundamental in information geometry
and in population genetics.

\begin{figure}[hbt]
  \centering \vspace{0.1in}
    \includegraphics[width=0.35\textwidth]{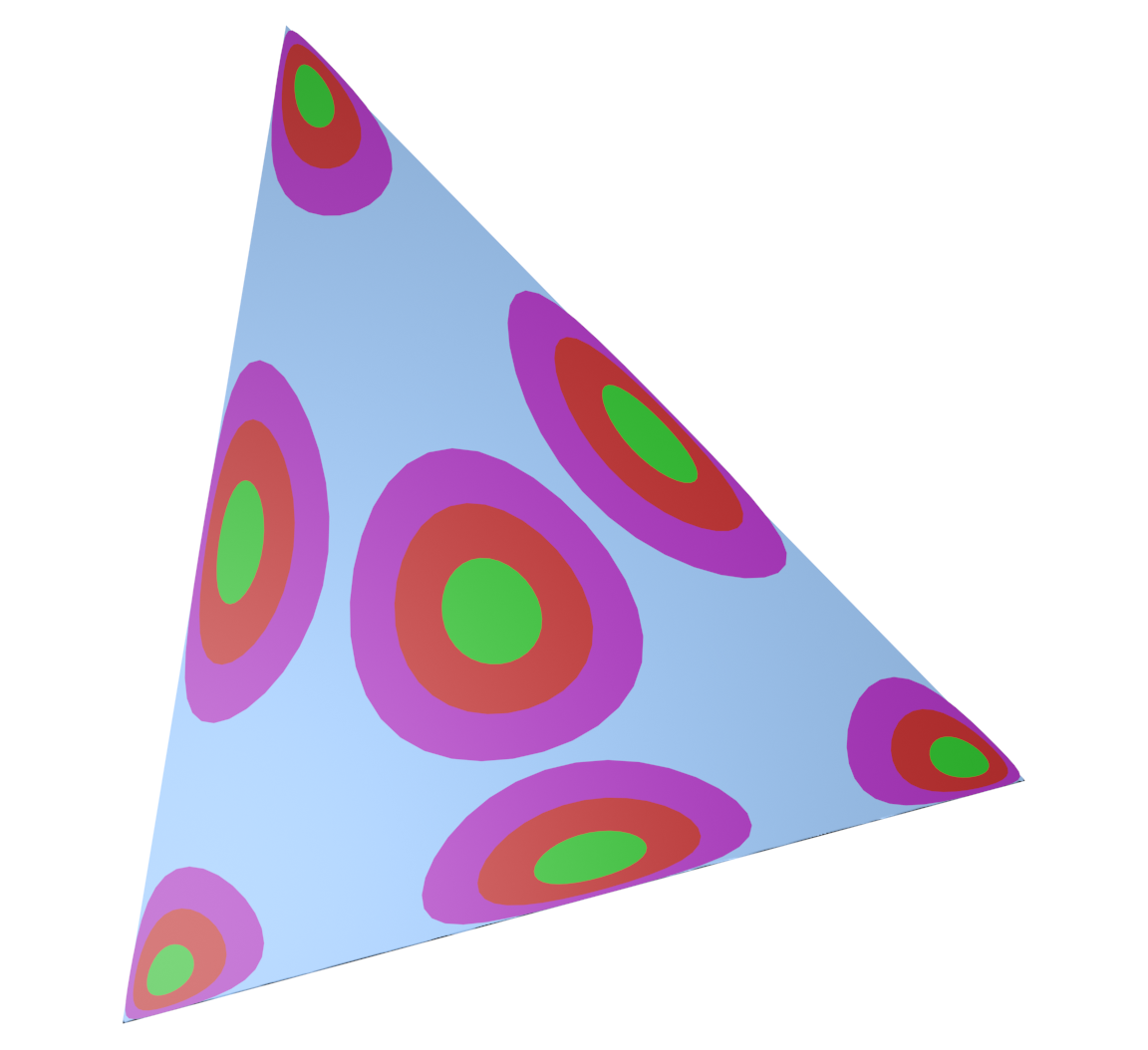}
    \hspace{0.25in}
    \includegraphics[width=0.35\textwidth]{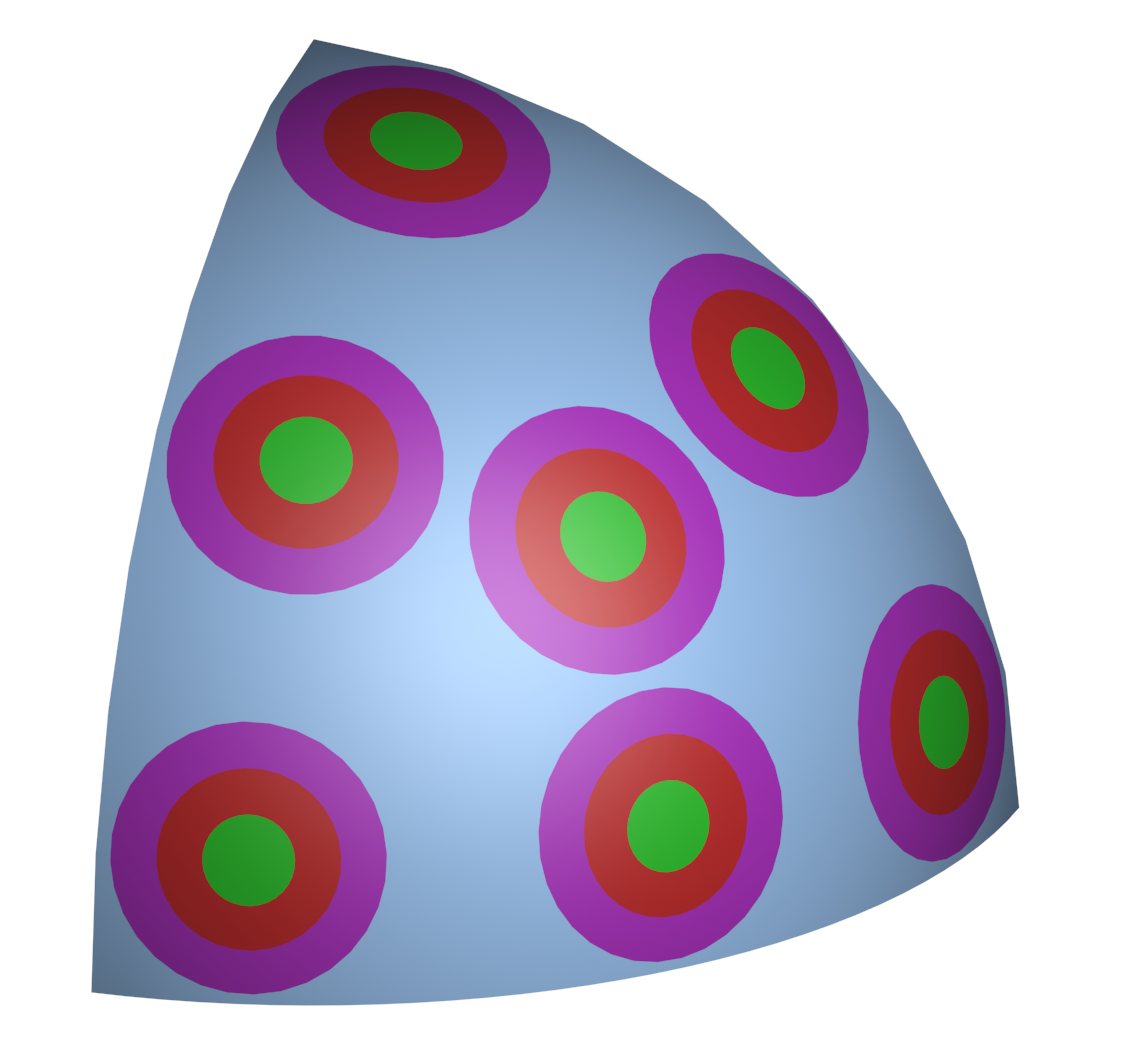}
    \vspace{0.1in}
  \caption{\emph{Left:} disk neighborhoods under the Fisher information metric
    of seven points in the standard triangle.
    \emph{Right:} the corresponding seven points and cap neighborhoods
    in the isometric non-negative octant of the $2$-sphere.
    For aesthetic reasons, the octant is scaled to $1 / \sqrt{2}$ times
    its actual size.}
  \label{fig:Fisher}
\end{figure}

To shed light on the Fisher information metric,
we map every point $\xxx = (x_0, x_1, \ldots, x_n)$ of $\Delta^n$
to the point $\varphi (\xxx) = (u_0, u_1, \ldots, u_n)$
with $u_i = \sqrt{2 x_i}$ for every $i$.
The coordinates of $\varphi (\xxx)$ are all non-negative
and satisfy $\sum_{i=0}^n u_i^2 = 2$.
In words, $\varphi (\xxx)$ is a point of $\sqrt{2} \Sspace_+^n$,
which is our notation for
the non-negative orthant of the sphere with radius $\sqrt{2}$ centered
at the origin in $\Rspace^{n+1}$;
see Figure \ref{fig:Fisher} on the right.
As noticed already by Antonelli \cite{Antonelli},
see also Akin \cite[page 39]{Aki79},
this mapping is an isometry between $\Delta^n$ and $\sqrt{2} \Sspace_+^n$.
We can therefore understand $\Delta^n$ under the Fisher information metric
by studying $\Sspace_+^n$ under the geodesic distance.
To get a handle on the difference between random sets in
$\Rspace^n$ and in $\Delta^n$, we compare
point sets selected from Poisson point processes
in $\Rspace^n$ and on $\Sspace^n$,
the latter being the topic of this article.
Figure \ref{fig:Fisher} illustrates the isometry by showing three
level lines each for seven points in the standard triangle on the left
and for the seven corresponding points in the positive orthant of the
  sphere on the right.

%%%%%%%%%%%%%%%%%%%%%%%%%%%%%%%%%%%%%%%%%%%%%%%%%%%%%%%%%%%%%%%%%%%%%%%%%%
\ourparagraph{Prior work.}
%%%%%%%%%%%%%%%%%%%%%%%%%%%%%%%%%%%%%%%%%%%%%%%%%%%%%%%%%%%%%%%%%%%%%%%%%%
This article builds on work from three different but related areas:
random polytopes, Poisson--Delaunay mosaics, and discrete Morse theory.

Consider the model in which a random polytope is generated by taking
the convex hull of randomly chosen points on the unit sphere.
The first paper with substantial results on this topic is Miles \cite{Mil71}.
The large body of work on the expected number of faces
of random polytopes and their volume is summarized and surveyed in \cite{BFV10,Hug13,Rei10,Sch08,ScWe08}.
A survey of recent results can be found in \cite{Ste14}.
The more general setting in which the points are selected on the boundary
of a convex body is addressed in \cite{ReSt16},
and the linear dependence of the expected number of faces on the number
of vertices is proved.

The study of Poisson--Delaunay mosaics in Euclidean space was started
by Miles with two seminal papers \cite{Mil70,Mil71} around 1970.
Considering the expected number of $k$-dimensional simplices
in an $n$-dimensional Poisson--Delaunay mosaic, he settles the question
for all values of $k$ in dimensions $n \leq 3$,
and for $k = n-1, n$ in any dimension $n$.
The first substantial extension of these results appeared in \cite{ENR16},
settling the question for all values of $k$ in dimension $n = 4$,
and determining the density of the radius of a typical simplex.
The main new idea in \cite{ENR16} is the classification of the simplices
based on the discrete Morse theory of the Delaunay mosaic,
and this approach is also central to the work in this paper.
Discrete Morse theory was first introduced as an abstract concept
in \cite{For98}, and its generalized version was used in \cite{BaEd16}
to study the radius function of a Delaunay mosaic.

%%%%%%%%%%%%%%%%%%%%%%%%%%%%%%%%%%%%%%%%%%%%%%%%%%%%%%%%%%%%%%%%%%%%%%%%%%
\ourparagraph{Concepts and notation.}
%%%%%%%%%%%%%%%%%%%%%%%%%%%%%%%%%%%%%%%%%%%%%%%%%%%%%%%%%%%%%%%%%%%%%%%%%%
Before stating our results, we introduce some concepts and notation;
the detailed description will follow in Section \ref{sec:3}.
We write $\Rspace^n$ for the $n$-dimensional \emph{Euclidean space},
$\Bspace^n \subseteq \Rspace^n$ for the closed \emph{unit ball},
and $\Sspace^{n-1} = \boundary{\Bspace^n}$ for the \emph{unit sphere},
the boundary of the unit ball.
Following \cite{ScWe08}, we write $\nu_n$ for the $n$-dimensional volume
of $\Bspace^n$ and $\sigma_n$ for the area (the $(n-1)$-dimensional volume)
of $\Sspace^{n-1}$.

To facilitate the comparison between Euclidean and spherical space,
we move up by one dimension and consider $\Sspace^n$,
which we equip with the \emph{geodesic distance},
$d \colon \Sspace^n \times \Sspace^n \to \Rspace$,
induced by the Euclidean metric in $\Rspace^{n+1}$;
see Section \ref{sec:3} for more details.
The relation between the geodesic distance and the Euclidean distance is
$\Gdist{x}{y} = 2 \arcsin \tfrac{\Edist{x}{y}}{2}$.
Let $\Capp{\GRad}{x} = \{ w \in \Sspace^n \mid \Gdist{w}{x} \leq \GRad \}$
be the \emph{spherical cap} with center $x \in \Sspace^n$
and \emph{geodesic radius} $0 \leq \GRad \leq \pi$. 
To measure the area of a cap, we use the \emph{Beta function},
$\Beta{a}{b} = \iBeta{1}{a}{b}$, and its \emph{incomplete} version,
$\iBeta{u}{a}{b} = \int_{t=0}^u t^{a-1} (1-t)^{b-1} \diff t$,
in which $0 \leq u \leq 1$.
For $\GRad \leq \tfrac{\pi}{2}$, the fraction of the sphere covered
by the cap is
$\Fraction{\GRad}  =  \tfrac{1}{2}
                      {\iBeta{s}{\tfrac{n}{2}}{\tfrac{1}{2}}} /
                      {\Beta{\tfrac{n}{2}}{\tfrac{1}{2}}}$,
in which $s = \sin^2 \GRad$ is the square of the Euclidean radius
measured in $\Rspace^{n+1}$; see \cite{Li11}.
The area of the cap is then
\begin{align}
  \Area{\GRad} &= \left\{ \begin{array}{rl}
                \Fraction{\GRad} \sigma_{n+1}
                  & \mbox{\rm for~} 0 \leq \GRad \leq \tfrac{\pi}{2}, \\
									\left[1 - \Fraction{\pi - \GRad}\right] \sigma_{n+1}
                  & \mbox{\rm for~} \tfrac{\pi}{2} \leq \GRad \leq \pi,
              \end{array} \right. 
  \label{eqn:CapArea}
\end{align}
in which $\Fraction{\pi - \GRad} = \Fraction{\GRad}$
because $\sin (\pi - \GRad) = \sin \GRad$.
Besides the Beta functions, we will use the \emph{Gamma function},
$\Gamma (k) = \gamma_\infty (k)$, and its \emph{lower incomplete} version,
$\iGama{u}{k} = \int_{t=0}^u t^{k-1} e^{-t} \diff t$,
in which $0 \leq u \leq \infty$.
The connection to the Beta functions is
$\Beta{a}{b} = \Gama{a} \Gama{b} / \Gama{a+b}$.
Finally, we write $\LGrass{k}{n}$ for the \emph{Grassmannian},
which consists of all $k$-dimensional planes that pass through the
origin in $\Rspace^n$.
This is a manifold of dimension $(n-k)\times k$ and measure
$\norm{\LGrass{k}{n}}
    = \tfrac{\sigma_n \cdot \sigma_{n-1} \cdot \ldots \cdot \sigma_{n-k+1}}
             {\sigma_1 \cdot \sigma_2 \cdot \ldots \cdot \sigma_k}$;
see \cite{ScWe08}.
This measure appears in the definition of constants that play an
important role in the statement of our results:
\begin{align}
  \Ccon{\ell}{k}{n}
    &=  \tfrac{\sigma_n \cdot \sigma_{n-1} \cdot \ldots \cdot \sigma_{n-k+1}}
              {\sigma_1 \cdot \sigma_2 \cdot \ldots \cdot \sigma_k}
        \cdot \tfrac{\Gama{k} n^{k-1} k!^{n-k} \sigma_k^{k+1}}
                    {(k+1) \sigma_n^k}
        \cdot \EAux{\ell}{k}{n} ,
  \label{eqn:constantC}                                          \\
	 \EAux{\ell}{k}{n} &= \Expected{\Volume{\uuu}^{n-k+1} \One{k-\ell}(\uuu)},
  \label{eqn:constantE}
\end{align}
in which $\uuu = (u_0, u_1, \ldots, u_k)$ is a sequence of $k+1$ points
chosen independently and uniformly at random on $\Sspace^{k-1}$,
with $\One{k-\ell} (\uuu) = 1$, if $k - \ell$ of the $k+1$ facets
of the $k$-simplex span $k$-planes that separate $0$ from $\uuu$,
and $\One{k-\ell} (\uuu) = 0$, otherwise; see \cite{ENR16},
where the same constants are studied in more detail.

\begin{figure}[hbt]
  \centering \vspace{0.1in}
  \includegraphics[width=0.5\textwidth]{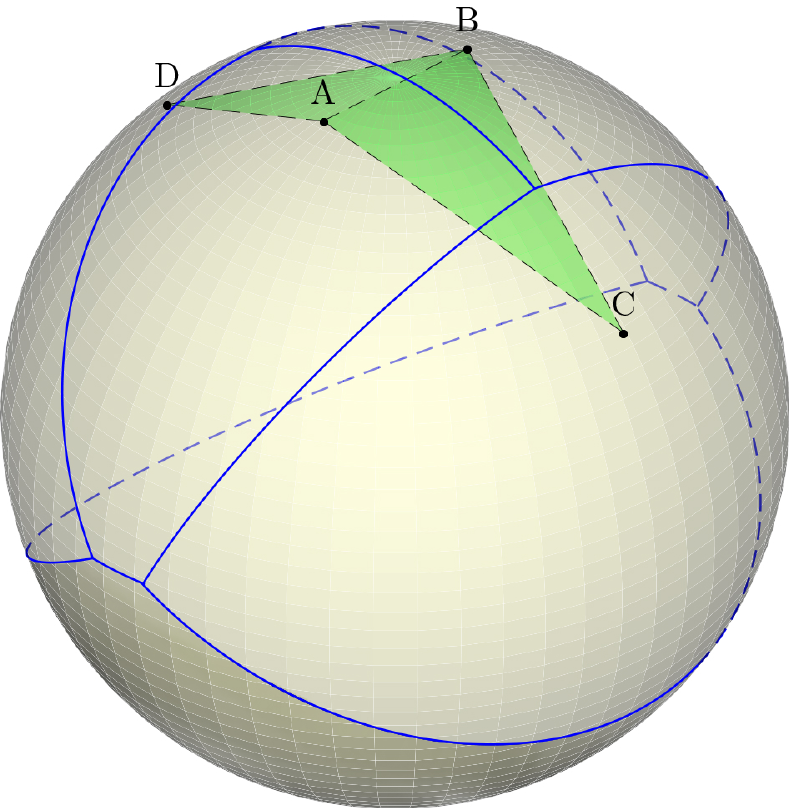} \vspace{0.1in}
  \caption{The Voronoi domains of four points on the $2$-dimensional sphere.
    The darker region in the south does not belong to any of these domains
    because the four points all belong to the northern hemisphere.
    The dual Delaunay complex consists of two triangles glued along a shared
    edge.}
  \label{fig:Voronoi}
\end{figure}

We follow \cite{Ren97} in defining the \emph{Voronoi domain} of a point
$x \in X$ as the set $\Voronoi{x}$ of points $w \in \Sspace^n$
that satisfy
$d(w,x) \leq d(w,y)$ for all $y \in X$ as well as $d(w,x) < \tfrac{\pi}{2}$.
Note that the Voronoi domains cover $\Sspace^n$ iff there is no closed
hemisphere that contains $X$; see Figure \ref{fig:Voronoi}.
The \emph{Delaunay mosaic} is the nerve of the Voronoi domains,
which we denote $\Delaunay{}{X}$.
Assuming general position, $\Delaunay{}{X}$ is isomorphic to a subcomplex
of the boundary complex of $\conv{X}$ in $\Rspace^{n+1}$,
and it is isomorphic to the entire boundary complex iff $X$
is not contained in any closed hemisphere.
Given a geodesic radius $0 \leq \GRad \leq \tfrac{\pi}{2}$,
we sometimes restrict the Voronoi domains to $\Voronoi{x} \cap \Capp{\GRad}{x}$,
for every $x \in X$,
and we write $\Delaunay{\GRad}{X}$ for the nerve of the thus
restricted Voronoi domains.
The \emph{(geodesic) radius function}, $\Rfun \colon \Delaunay{}{X} \to \Rspace$,
maps every Delaunay simplex to the radius of its
\emph{smallest empty circumscribed cap};
see details in Section \ref{sec:3}.
By assumption, this radius is always less than $\tfrac{\pi}{2}$,
and we observe that $\Rfun^{-1} [0, \GRad] = \Delaunay{\GRad}{X}$
for all $\GRad \geq 0$.
For points in general position, $\Rfun$ is a generalized discrete
Morse function, as defined in \cite{BaEd16,For98},
which implies a partition of $\Delaunay{}{X}$ into maximal \emph{intervals}
$[L,U] = \{Q \mid L \subseteq Q \subseteq U\}$
consisting of simplices with equal function value, as discussed in Section \ref{sec:3}.
The \emph{type} of the interval $[L, U]$ is $(\ell, k)$,
in which $\ell = \dime{L}$ and $k = \dime{U}$.
If $L = U$, then the interval contains a single simplex,
which we call a \emph{critical simplex} of $\Rfun$.

We need one more concept to express the asymptotic behavior of the
expected numbers, when their density goes to infinity.
Assuming a Poisson point process with density $\density > 0$ on $\Sspace^n$,
for a cap with geodesic radius $\GRad$, we call
$\NRad = \GRad \density^{1/n}$ the \emph{normalized radius} of the cap.
It is the geodesic radius of the cap after scaling the unit sphere
to the sphere with area $\density \sigma_n$.

%%%%%%%%%%%%%%%%%%%%%%%%%%%%%%%%%%%%%%%%%%%%%%%%%%%%%%%%%%%%%%%%%%%%%%%%%%
\ourparagraph{Results.}
%%%%%%%%%%%%%%%%%%%%%%%%%%%%%%%%%%%%%%%%%%%%%%%%%%%%%%%%%%%%%%%%%%%%%%%%%%
Our main result concerns the radius function on the Delaunay mosaic.
For each type, we express the expected number of intervals of that type with
normalized radius smaller than a threshold in terms of an integral,
which we evaluate asymptotically, when the density of the Poisson point
process goes to infinity.
\begin{theorem}[Main Result]
  \label{thm:MainResult}
  Let $X$ be a Poisson point process with density $\density > 0$ on $\Sspace^n$.
  For any integers $1 \leq \ell \leq k \leq n$ and any real number
  $0 < \GRad_0 < \tfrac{\pi}{2}$, the expected number of intervals
  of type $(\ell, k)$ and geodesic radius at most $\GRad_0$ is
  \begin{align}
    \Expected{\ccon{\ell}{k}{n}, \GRad_0} &= 
      \density \sigma_{n+1} \cdot \tfrac{\sigma_n^k}
                                        {2 \Gamma(k) n^{k-1}}
      \cdot \Ccon{\ell}{k}{n}
      \int\displaylimits_{t = 0}^{s}
        \density^k t^{\frac{kn-2}{2}} (1-t)^{\frac{n-k-1}{2}}
        \Pempty{\sqrt{t}} \diff t,
    \label{eqn:MainResult}
  \end{align}
  in which $s = \sin^2 \GRad_0$ is the square of the maximum Euclidean radius, 
  and $\Pempty{\ERad}$ is the probability that a spherical cap
  with geodesic radius $\GRad = \arcsin \ERad$ contains no points of $X$,
  namely $\Pempty{\ERad} = e^{- \density \Area{\GRad}}$.
  Let now $\density \to \infty$.
  For any $\NRad_0 \in [0, +\infty]$,
  the expected number of intervals of type $(\ell, k)$
  and normalized radius at most $\NRad_0$ is
  \begin{align}
    \Expected{\ccon{\ell}{k}{n}, \NRad_0} 
      &=  \density \sigma_{n+1} \cdot \tfrac{\iGama{v}{k}}{\Gama{k}}
	    \cdot \Ccon{\ell}{k}{n} + \littleoh{\density},
    \label{eqn:asymp_res}
  \end{align}
  in which $v = \NRad_0^n \nu_n$ is the volume of the $n$-ball with
  radius $\NRad_0$.
\end{theorem}

\noindent {\sc Remarks.}
({\sc 1a}) Theorem \ref{thm:MainResult}
does not cover the case $\ell = 0$, i.e., intervals containing vertices,
but here the results are straightforward.
Specifically, the expected number of critical vertices is
$\Expected{\ccon{0}{0}{n}, \GRad_0} = \density \sigma_{n+1}$,
for every $\GRad_0 \geq 0$,
and $\ccon{0}{k}{n} = 0$ for every $k \geq 1$.

({\sc 1b}) We will prove that for constant $s$,
the integral in \eqref{eqn:MainResult} is bounded away from
both $0$ and $\infty$.
This implies that the expected number of intervals
in \eqref{eqn:MainResult} is of order $\bigtheta{\density}$;
compare with \cite{ReSt16}.

({\sc 1c}) We will also prove that setting $\NRad_0 = \infty$ in \eqref{eqn:asymp_res}
gives the total number of intervals of type $(\ell, k)$ as
$\Expected{\ccon{\ell}{k}{n}}
    =  \density \sigma_{n+1} \cdot \Ccon{\ell}{k}{n} + \littleoh{\density}$.
On the other hand, letting $\NRad_0 \to \infty$,
we get the total number of intervals of geodesic radius $\bigtheta{\density^{- 1/n}}$.
This implies that the number of intervals with
radius $\littleomega{\density^{- 1/n}}$ is $\littleoh{\density}$.
Note that also the number of intervals with radius
$\littleoh{\density^{- 1/n}}$ is $\littleoh{\density}$.

\medskip
The total number of simplices of dimension $j$ in the Delaunay mosaic
is easy to deduce from the number of intervals:
$\dcon{j}{n} = \sum_{k=j}^n \sum_{\ell=0}^j \binom{k-\ell}{k-j} \ccon{\ell}{k}{n}$.
Accordingly, we define the constant
$\Dcon{j}{n} = \sum_{k=j}^n \sum_{\ell=0}^j \binom{k-\ell}{k-j} \Ccon{\ell}{k}{n}$.
We generalize this relation so it depends on a normalized radius threshold:
\begin{corollary}[Delaunay Simplices]
  \label{cor:DelaunaySimplices}
  Let $X$ be a Poisson point process with density $\density > 0$ on $\Sspace^n$.
  For any integer $j \geq 1$ and any non-negative real number $\NRad_0$,
  the expected number of $j$-simplices of $\Delaunay{}{X}$ with
  normalized radius at most $\NRad_0$ is
  \begin{align}
    \Expected{\dcon{j}{n}, \NRad_0}  &=
      \density \sigma_{n+1} \cdot \sum_{k=j}^n \tfrac{\iGama{v}{k}}{\Gama{k}}
      \sum_{\ell=0}^j \binom{k-\ell}{k-j} \Ccon{\ell}{k}{n}
          + \littleoh{\density},
    \label{eqn:delaunay_simplices}
  \end{align}
  in which $v = \NRad_0^n \nu_n$.
  Setting
  \begin{align}
    G_j^n (\NRad_0) 
      &=  \sum_{k=j}^n \tfrac{\iGama{v}{k}}{\Gama{k}}
          \sum_{\ell=0}^j \binom{k-\ell}{k-j} \tfrac{\Ccon{\ell}{k}{n}}{\Dcon{j}{n}},
    \label{eqn:distribution}
  \end{align}
  we thus get the distribution of the normalized radius of the typical $j$-simplex
  in the limit when $\density \to \infty$.
\end{corollary}

\noindent {\sc Remarks.}
({\sc 2a}) Observe that $\density \sigma_{n+1}$ is the expected number
of points in $X$.
Comparing with \cite{ENR16}, we thus notice that
\eqref{eqn:asymp_res}, \eqref{eqn:delaunay_simplices},
and \eqref{eqn:distribution}
are essentially the same expressions as for the Poisson point
process in $\Rspace^n$.
This justifies the title of this article.

({\sc 2b}) While we state our results for Poisson point processes,
very similar results can be obtained for the uniform distribution;
see Appendix \ref{app:uniform}.  

%%%%%%%%%%%%%%%%%%%%%%%%%%%%%%%%%%%%%%%%%%%%%%%%%%%%%%%%%%%%%%%%%%%%%%%%
\ourparagraph{Outline.}
%%%%%%%%%%%%%%%%%%%%%%%%%%%%%%%%%%%%%%%%%%%%%%%%%%%%%%%%%%%%%%%%%%%%%%%%
Section \ref{sec:2} introduces the main technical tool used to prove
our results.
Section \ref{sec:3} gives the background, including discrete Morse theory.
Section \ref{sec:4} proves the integral equation and the asymptotic result
both stated in Theorem \ref{thm:MainResult}.
Section \ref{sec:5} concludes the paper.

%\newpage
%%%%%%%%%%%%%%%%%%%%%%%%%%%%%%%%%%%%%%%%%%%%%%%%%%%%%%%%%%%%%%%%%%%%%%%%%%
%%%%%%%%%%%%%%%%%%%%%%%%%%%%%%%%%%%%%%%%%%%%%%%%%%%%%%%%%%%%%%%%%%%%%%%%%%
\section{Blaschke{\textendash}Petkantschin Formula for the Sphere}
\label{sec:2}
%%%%%%%%%%%%%%%%%%%%%%%%%%%%%%%%%%%%%%%%%%%%%%%%%%%%%%%%%%%%%%%%%%%%%%%%%%
%%%%%%%%%%%%%%%%%%%%%%%%%%%%%%%%%%%%%%%%%%%%%%%%%%%%%%%%%%%%%%%%%%%%%%%%%%

This section introduces a formula of Blaschke--Petkantschin type
used in the proof of Theorem \ref{thm:MainResult}.
Since it is a stand-alone result, not specific to the problem addressed
in this article, we present it before discussing the background related
to the subject in this paper.
In its basic form, the \emph{Blaschke--Petkantschin formula} writes
an integral over $\Rspace^{n+1}$ as an integral over the Grassmannian,
$\LGrass{k}{n+1}$.
We adapt the original such formula to the $n$-sphere.
Formulas of this type were studied in \cite{Za90}.
To express the result, we write $\plane^\perp$ for the $(n-k+1)$-plane
orthogonal to the $k$-plane $\plane$,
both passing through the origin in $\Rspace^{n+1}$,
and we write $S_\plane$ for the unit $(k-1)$-sphere in $\plane$.
As usual, we use boldface to denote sequences of points:
$\xxx = (x_0, x_1, \ldots, x_k)$, etc.
A shortcut $p+r\uuu$ is used for $(p+r u_0, p+ru_1, \ldots, p+ru_k)$.
The integrations are with respect to the standard
measure on the Grassmanian and the Lebesgue measures in the plane and on the sphere.
\begin{theorem}[Blaschke--Petkantschin for the Sphere]
  \label{thm:BPforSphere}
  Let $n$ be a positive integer, $1 \leq k \leq n$,
  and $f \colon (\Sspace^n)^{k+1} \to \Rspace$ a non-negative measurable function.
  Then
  \begin{align}
    \!\! \int\displaylimits_{\xxx \in (\Sspace^n)^{k+1}}
          \!\!\! f(\xxx) \diff \xxx
      &=  \int\displaylimits_{\plane \in \LGrass{k}{n+1}}
          \int\displaylimits_{p \in \plane^\perp} \!\! \ERad^{kn-2}
          \!\! \int\displaylimits_{\uuu \in (S_\plane)^{k+1}}
            \!\!\! f(p + \ERad \uuu) \left[ k!\Volume{\uuu} \right]^{n-k+1}
          \diff \uuu \diff p \diff \plane ,
    \label{eqn:BPforSphere}
  \end{align}
  in which $\ERad^2 = 1-\norm{p}^2$, implicitly assuming $\norm{p} \leq 1$,
  and $\Volume{\uuu}$ is the $k$-dimensional volume of the convex hull
  of the points in $\uuu$, which is a $k$-simplex.
  If $f$ is rotationally symmetric, we define
  $f_\ERad (\uuu) = f(p + \ERad \uuu)$, in which
  $\uuu$ is a $k$-simplex on $\Sspace^{k-1} \subseteq \Rspace^k$,
  and $p$ is any point with $\norm{p}^2 = 1 - \ERad^2 \leq 1$ in the $(n-k+1)$-plane
  orthogonal to $\Rspace^k \subseteq \Rspace^{n+1}$.
  With this notation, we have
  \begin{align}
    \!\! \int\displaylimits_{\xxx \in (\Sspace^n)^{k+1}}
          \!\!\!\!\! f(\xxx) \diff \xxx
      &=  \tfrac{\sigma_{n+1}}{2} \norm{\LGrass{k}{n}}
          \int\displaylimits_{t=0}^1
            t^{\frac{kn-2}{2}}\!(1\! -\! t)^{\frac{n-k-1}{2}}
          \!\!\!\!\!\!\! \int\displaylimits_{\uuu \in (\Sspace^{k-1})^{k+1}}
            \!\!\!\!\!\!\!\! f_{\sqrt{t}} (\uuu) \left[ k!\Volume{\uuu} \right]^{n-k+1}
          \! \diff \uuu \diff t .
    \label{eqn:rotationally-symmetric}
  \end{align}
\end{theorem}
\ourproof
  We first argue that $f$ may be assumed to be continuous.
  Consider the subset $M$ of
  $\LGrass{k}{n+1} \times \Rspace^{n+1} \times (\Rspace^{n+1})^{k+1}$
  consisting of all triplets $(P, p, \uuu)$
  such that $p \in P^\perp$, $\norm{p} < 1$, and $\uuu \in (S_P)^{k+1}$.
  Clearly, $M$ is a submanifold of the product space with a natural measure.
  Recall that $\ERad^2 = 1-\norm{p}^2$ and
  consider the mapping $T \colon M \to (\Sspace^n)^{k+1}$
  defined by $T(P, p, \uuu) = p + \ERad \uuu$.
  It is a bijection up to a set of measure $0$.
  By Theorem 20.3 in \cite{HeSt65},
  there exists a corresponding Jacobian $J \colon M \to \Rspace$,
  meaning that every integrable function $f$ satisfies
  $\int\displaylimits_{\xxx \in (\Sspace^n)^{k+1}} f(\xxx) \diff \xxx
    = \int\displaylimits_{y \in M} f(T(y)) J(y) \diff y$.
  For non-negative $f$, the right-hand side integral can be split
  using Fubini's theorem.
  The existence of the Jacobian is thus settled,
  and to find its values, we may assume that $f$ be continuous.
  
\begin{figure}[hbt]
    \centering \vspace{0.1in}
    \resizebox{!}{2.0in}{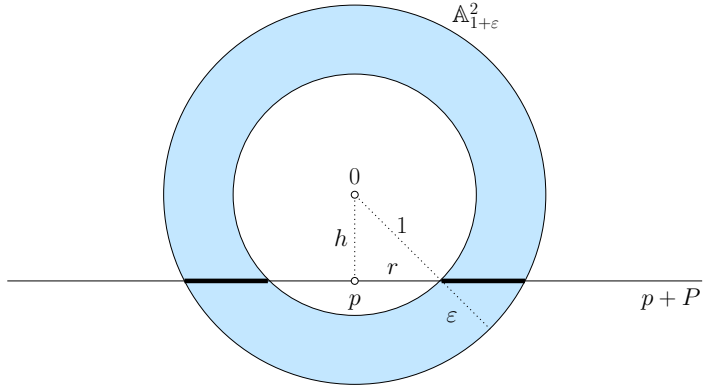} \vspace{0.1in}
    \caption{For $h = \norm{p} < 1$, the slice of the $(n+1)$-dimensional
      annulus is a $k$-dimensional annulus.  In this picture,
      $n+1 = 2$ and $k = 1$.}
    \label{fig:annulus}
  \end{figure}

  The main idea in the rest of the proof
  is to thicken $\Sspace^n$ to an $(n+1)$-dimensional annulus,
  to apply the original Blaschke--Petkantschin formula to this annulus,
  and to take the limit when we shrink the annulus back to $\Sspace^n$.
  We write $\Aspace^{n+1}_{1+\ee} =
            (1+\ee) \Bspace^{n+1} \setminus \interior{\Bspace^{n+1}}$
  for the $(n+1)$-dimensional annulus with inner radius $1$ and
  outer radius $1+\ee$.
  We begin by extending $f$ from the sphere to the annulus.
  Specifically, for points $y_i \in \Aspace^{n+1}_{1+\ee}$, we set
  \begin{align}
    F(y_0,y_1,\ldots,y_k)  &=  f \left( {y_0}/{\norm{y_0}},
                                        {y_1}/{\norm{y_1}}, \ldots,
                                        {y_k}/{\norm{y_k}} \right) .
  \end{align}
  Since $f$ is continuous on the $(k+1)$-fold product of spheres,
  by assumption, $F$ is continuous on the $(k+1)$-fold product of annuli.
  Because $F$ is continuous on a compact set and therefore bounded
  and uniformly continuous, we have
  \begin{align}
    \int\displaylimits_{\xxx \in (\Sspace^n)^{k+1}}
          f(\xxx) \diff \xxx
      &=  \lim_{\ee \to 0} \tfrac{1}{\ee^{k+1}}
             \int\displaylimits_{\yyy \in (\Aspace^{n+1}_{1+\ee})^{k+1}}
             F(\yyy) \diff \yyy
    \label{eqn:innermost1} \\
      &=  \lim_{\ee \to 0} \tfrac{1}{\ee^{k+1}}
            \int\displaylimits_{\plane \in \LGrass{k}{n+1}}
            \int\displaylimits_{p \in \plane^\perp}
            \int\displaylimits_{\uuu \in A^{k+1}}
              F(\uuu) [k! \Volume{\uuu}]^{n-k+1}
              \diff \uuu \diff p \diff \plane ,
    \label{eqn:innermost2}
  \end{align}
  in which $A = \Aspace^{n+1}_{1+\ee} \cap [p+\plane]$
  is the $k$-dimensional slice of the $(n+1)$-dimensional annulus
  defined by $\plane$ and $p$.
  We obtain \eqref{eqn:innermost2} from \eqref{eqn:innermost1} by applying
  the standard Blaschke--Petkantschin formula in $\Rspace^{n+1}$
  to the function $F(\yyy)$ times the indicator function of
  the $(k+1)$-fold product of annuli,
  and then absorb the indicator into the integration domain.
  To continue, we investigate the slice of the annulus whose
  $(k+1)$-fold product is the innermost integration domain;
  see Figure \ref{fig:annulus}.
  Write $h = \norm{p}$ for the height of the slice, which is non-empty
  for $0 \leq h \leq 1+\ee$. $A$ is a (possibly degenerate) $k$-dimensional
  annulus, with squared inner radius $\ERad^2 = \max \{0, 1-h^2\}$
  and squared outer radius $\ERad_\ee^2 = (1+\ee)^2-h^2$.
  We split the integration domain into three regions:
  $h \leq 1-\ee^{-0.2}$, $1 - \ee^{0.2} < h \leq 1$, and $1 < h \leq 1 + \ee$.
  
  We first show that the contribution of the region
  $1 - \ee^{0.2} < h \leq 1$ is small.
  To get started, note that
  $\ERad_\ee - \ERad  =  (\ERad_\ee^2 - \ERad^2) / (\ERad_\ee + \ERad)
                      =  (2\ee + \ee^2) / (\ERad_\ee + \ERad)$.
  For small $\ee$, this implies
  $\ERad_\ee - \ERad \leq \const \cdot \ee / \ERad_\ee$,
  in which we deliberately avoid the computation of the constant.
  With this, we can bound the $k$-dimensional volume of $A$.
  Assuming $k \geq 2$, we get
  $\Volume{A}  =  \nu_k (\ERad_\ee^k - \ERad^k)                            
               =  \nu_k (\ERad_\ee - \ERad)
                    (\ERad_\ee^{k-1} + \ERad_\ee^{k-2}\ERad + \ldots + \ERad^{k-1})
             \leq  \const \cdot \ee \ERad_\ee^{k-2}$,
  in which the constant depends only on $k$ and $n$.
  As noted before, the inequality also holds for $k = 1$.
  Since $h > 1 - \ee^{0.2}$, we also get
  $\ERad_\ee^2  <  (1+\ee)^2 - (1-\ee^{0.2})^2
             \leq  \ee^2 + 2 \ee + 2 \ee^{0.2} - \ee^{0.4}$ for small $\ee$,
  which implies $\ERad_\ee < \const \cdot \ee^{0.1}$.
  Clearly, the $k$-dimensional volume of any $k$-simplex
  with vertices inside $A$ can not exceed a constant times the $k$-th power
  of the diameter of $A$, which is $2\ERad_\ee$,
  implying $\Volume{\uuu} \leq \const \cdot \ERad_\ee^k$.
  Recalling that $F$ is bounded, we thus get
  \begingroup
    \allowdisplaybreaks
    \begin{align}
      \MoveEqLeft \left|\,\int\displaylimits_{\plane \in \LGrass{k}{n+1}}\,\,
          \int\displaylimits_{\substack{p \in \plane^\perp \\
                              \|p\|<1-\ee^{0.2}}} \tfrac{1}{\ee^{k+1}}
        \int\displaylimits_{u\in A^{k+1}} F(\uuu) [k! \Volume{\uuu}]^{n-k+1}
          \diff \uuu \diff p \diff\plane \right|                             \\
      &\leq \const \int\displaylimits_{h=1-\ee^{0.2}}^1
        \tfrac{1}{\ee^{k+1}} \Volume{A}^{k+1} \Volume{\uuu}^{n-k+1} \diff h  \\
      &\leq \const \int\displaylimits_{h=1-\ee^{0.2}}^1 \tfrac{1}{\ee^{k+1}}
        (\ee \ERad_\ee^{k-2})^{k+1} \ERad_\ee^{k(n-k+1)} \diff h             \\
      &\leq \const \int\displaylimits_{h=1-\ee^{0.2}}^1 \ERad_\ee^{kn-2} \diff h
        \;\leq\; \const \cdot \ee^{0.2} \cdot \ee^{0.1 (kn-2)} \;\to\; 0.
    \label{eqn:cutoff_one}
  \end{align}
  \endgroup
  Here we use the bound on $\ERad_\ee$ for the last inequality,
  and $kn \geq 1$ to see that the expression tends to zero.
  Next consider the region $1 < h \leq 1+\ee$,
  in which $A$ is a ball of radius $\ERad_\ee$,
  so $\Volume{A} = \nu_k \ERad_\ee^k$.
  We have $\Volume{\uuu} \leq \nu_k \ERad_\ee^k$, as before,
  and $\ERad_\ee^2 \leq (1+\ee)^2 - 1$,
  which implies $\ERad_\ee \leq \const \cdot \sqrt{\ee}$.
  With this, we can again establish the vanishing of the integral as $\ee \to 0$:
  \begin{align}
    \MoveEqLeft \left|\,\int\displaylimits_{\plane \in \LGrass{k}{n+1}}\,\,
      \int\displaylimits_{\substack{p \in \plane^\perp \\
                                    1 \leq \norm{p} \leq 1+\ee}} \tfrac{1}{\ee^{k+1}}
        \int\displaylimits_{u\in A^{k+1}} F(\uuu) [k! \Volume{\uuu}]^{n-k+1}
          \diff \uuu \diff p\diff\plane \right|                              \\
      &\leq \const \int\displaylimits_{h=1}^{1+\ee}
        \tfrac{1}{\ee^{k+1}} \Volume{A}^{k+1} \Volume{\uuu}^{n-k+1} \diff h  \\
      &\leq \const \int\displaylimits_{h=1}^{1+\ee} \tfrac{1}{\ee^{k+1}}
        \ERad_\ee^{k(n+2)} \diff h
      \;\leq\; \const \cdot \ee \cdot \ee^{(kn-2)/2} \;\to\; 0  .
    \label{eqn:cutoff_two}
  \end{align}
  We have thus established that the relevant region is
  $0 \leq h \leq 1-\ee^{0.2}$, and we are ready to investigate its contribution.
  First, we claim that the width of the annulus $A$ is
  \begin{align}
    \ERad_\ee - \ERad  &=  \ERad \sqrt{1 + \tfrac{2\ee+\ee^2}{\ERad^2}} - \ERad
                =  \tfrac{\ee}{\ERad} + \littleoh{\ee} .
    \label{eqn:width}
  \end{align}
  To get the right-hand side of \eqref{eqn:width}, we use the Taylor
  expansion of
  $g(x) = (1+x)^{1/2} = 1 + \tfrac{1}{2} x - \tfrac{1}{2} x^2 + \ldots$,
  and $\ERad > \ee^{0.1}$ as well as $x = (2 \ee + \ee^2)/\ERad^2 < 3 \ee^{0.8}$,
  which we get from the assumed $h \leq 1-\ee^{0.2}$.
  observing that $\ee^2/(2 \ERad^2) = \bigoh{\ee^{1.8}}$, we get
  $\ERad g(x) - \ERad
    = \tfrac{\ee}{\ERad} + \bigoh{\ERad \ee^{1.8}} + \bigoh{\ERad \ee^{1.6}}$
  and therefore \eqref{eqn:width}.
  Using the fact that $F(\uuu)$ is equal to $f(\uuu)$
  when all points lie on the inner sphere and the uniform continuity of $F$ and
  writing $S_r$ for the $(k-1)$-sphere with center $p$ and radius $r$
  in $\plane \in \LGrass{k}{n+1}$, we get
  \begin{align}
    \int\displaylimits_{\uuu \in A^{k+1}}
          \!\! \tfrac{1}{\ee^{k+1}}F(\uuu) [k! \Volume{\uuu}]^{n-k+1} \diff \uuu
      &=  \left( \tfrac{1}{\ERad} \right)^{k+1}
          \!\!\!\! \int\displaylimits_{\uuu \in (S_{r})}^{k+1}
          \!\! f(\uuu) [k! \Volume{\uuu}]^{n-k+1} \diff \uuu
          + \littleoh{1},
    \label{eqn:flatannulus}
  \end{align}
  in which the integration domain on the right is the
  $k$-fold product of the $(k-1)$-sphere
  with center $p$ and radius $\ERad$ in $\plane$,
  and $\littleoh{1}$ is uniform over $p$ and $\plane$.
  Substituting \eqref{eqn:cutoff_one}, \eqref{eqn:cutoff_two},
  and \eqref{eqn:flatannulus} into \eqref{eqn:innermost2}, we finally get
  \begin{align}
    \MoveEqLeft \int\displaylimits_{\xxx \in (\Sspace^n)^{k+1}} f(\xxx) \diff \xxx \\
      &=  \lim_{\ee \to 0}
          \int\displaylimits_{\plane \in \LGrass{k}{n+1}}
          \!\!\!\int\displaylimits_{\substack{p \in \plane^\perp \\
                                    \|p\| \leq 1-\ee^{0.2}}}\!\!
            \!\left[ \tfrac{1}{\ERad^{k+1}}
          \!\!\!\!\!\!\int\displaylimits_{\uuu \in (S_{r})^{k+1}}\!\!\!\!
             \!\!\! f(\uuu) [k! \Volume{\uuu}]^{n-k+1} \diff \uuu + \littleoh{1} \right]
             \! \diff p \diff \plane + \littleoh{1} 
      \label{eqn:BP0} \\
      &=  \int\displaylimits_{\plane \in \LGrass{k}{n+1}}
          \int\displaylimits_{p \in \plane^\perp} \left( \tfrac{1}{\ERad} \right)^{k+1}
          \!\!\! \int\displaylimits_{\uuu \in (S_{r})^{k+1}}
             \!\!\! f(\uuu) [k! \Volume{\uuu}]^{n-k+1} \diff \uuu \diff p \diff \plane 
      \label{eqn:BP1} \\
      &=  \int\displaylimits_{\plane \in \LGrass{k}{n+1}}
          \int\displaylimits_{p \in \plane^\perp} \ERad^{kn-2}
          \!\!\! \int\displaylimits_{\uuu \in (S_\plane)^{k+1}}
             \!\!\! f(p + \ERad \uuu) [k! \Volume{\uuu}]^{n-k+1}
             \diff \uuu \diff p \diff \plane ,
      \label{eqn:BP2}
  \end{align}
  in which we drop the $\norm{p} \leq 1-\ee^{0.2}$ condition in \eqref{eqn:BP0}
  for the implicitly assumed $\norm{p} \leq 1$
  when passing to \eqref{eqn:BP1},
  which we can do because the difference vanishes in the limit
  and \eqref{eqn:BP2} is obtained by rescaling and
  translating the sphere in \eqref{eqn:BP1}.
  Indeed, the power of $\ERad$ is a consequence of scaling the volume of the
  $k$-simplex, adjusting the volume of the integration domain,
  and subtracting the power we have already in \eqref{eqn:BP1}:
  $k(n-k+1)+(k-1)(k+1)-(k+1) = kn-2$.
  This proves the first relation claimed in Theorem \ref{thm:BPforSphere}.

  To get the second relation, we simplify the first by exploiting the
  rotational symmetry of $f$.
  Recalling that $\ERad^2 = 1 - \norm{p}^2$, it makes sense to define
  $f_\ERad (\uuu) = f(p + \ERad \uuu)$ on the $(k+1)$-fold product of
  $S_\plane \subseteq \Sspace^n$
  because the direction of $p$ does not matter for a fixed height.
  Neither does $\plane$ influence the function for a fixed height,
  so we can define $f_\ERad$ on $(\Sspace^{k-1})^{k+1}$.
  Thus
  \begin{align}
    \!\!\!\! \int\displaylimits_{\xxx \in (\Sspace^n)^{k+1}}
          \!\!\! f(\xxx) \diff \xxx
      &=  \norm{\LGrass{k}{n+1}}
          \!\! \int\displaylimits_{p \in \Bspace^{n-k+1}}
            \!\! \ERad^{kn-2}
          \!\! \int\displaylimits_{\uuu \in (\Sspace^{k-1})^{k+1}}
            \!\! f_\ERad (\uuu) [k! \Volume{\uuu}]^{n-k+1} \diff \uuu \diff p
      \label{eqn:BPsymm1} \\
      &=  \norm{\LGrass{k}{n+1}} \sigma_{n-k+1}
          \!\! \int\displaylimits_{h=0}^1
            \!\! h^{n-k} \ERad^{kn-2}
          \!\! \int\displaylimits_{\uuu \in (\Sspace^{k-1})^{k+1}}
            \!\! f_\ERad (\uuu) [k! \Volume{\uuu}]^{n-k+1} \diff \uuu \diff h
      \label{eqn:BPsymm2} \\
      &=  \tfrac{\sigma_{n+1}}{2} \norm{\LGrass{k}{n}}
          \! \int\displaylimits_{t=0}^1 t^{\frac{kn-2}{2}} (1-t)^{\frac{n-k-1}{2}}
          \!\!\!\!\!\! \int\displaylimits_{\uuu \in (\Sspace^{k-1})^{k+1}}
            \!\!\!\!\!\!\!\!\! f_\ERad (\uuu) [k! \Volume{\uuu}]^{n-k+1} \diff \uuu \diff t ,
      \label{eqn:BPsymm3}
  \end{align}
  in which $t = \ERad^2 = 1-h^2$.
  We get \eqref{eqn:BPsymm1} from \eqref{eqn:BP2} because
  every $\plane \in \LGrass{k}{n+1}$ contributes the same to the integral.
  Similarly, we get \eqref{eqn:BPsymm2} from \eqref{eqn:BPsymm1}
  by integrating over the range of heights and compensating for the different
  sizes of the corresponding spheres,
  aka expressing the integral in polar coordinates.
  Finally, we get \eqref{eqn:BPsymm3} from \eqref{eqn:BPsymm2} by
  substituting $t$ for $\ERad^2$, $1-t$ for $h^2$, and $\diff t$ for $-2 h\diff h$,
  noting that the minus sign is absorbed by reversing the limits of integration.
  This proves the second relation in Theorem \ref{thm:BPforSphere}.
\eop

%\newpage
%%%%%%%%%%%%%%%%%%%%%%%%%%%%%%%%%%%%%%%%%%%%%%%%%%%%%%%%%%%%%%%%%%%%%%%%%%
%%%%%%%%%%%%%%%%%%%%%%%%%%%%%%%%%%%%%%%%%%%%%%%%%%%%%%%%%%%%%%%%%%%%%%%%%%
\section{Background}
\label{sec:3}
%%%%%%%%%%%%%%%%%%%%%%%%%%%%%%%%%%%%%%%%%%%%%%%%%%%%%%%%%%%%%%%%%%%%%%%%%%
%%%%%%%%%%%%%%%%%%%%%%%%%%%%%%%%%%%%%%%%%%%%%%%%%%%%%%%%%%%%%%%%%%%%%%%%%%

This section introduces the geometric background needed to appreciate
the results in this paper.
After presenting the diagrams under study,
we explain the connection to discrete Morse theory,
and finally describe how we generate random diagrams.

%%%%%%%%%%%%%%%%%%%%%%%%%%%%%%%%%%%%%%%%%%%%%%%%%%%%%%%%%%%%%%%%%%%%%%%%%%
\ourparagraph{Voronoi tessellations and Delaunay mosaics.}
%%%%%%%%%%%%%%%%%%%%%%%%%%%%%%%%%%%%%%%%%%%%%%%%%%%%%%%%%%%%%%%%%%%%%%%%%%
We recall that the object under consideration is
$\Sspace^n \subseteq \Rspace^{n+1}$ with the \emph{geodesic distance},
$\GG \colon \Sspace^n \times \Sspace^n \to \Rspace$,
the metric inherited from the Euclidean metric on $\Rspace^{n+1}$.
The distance between any pair of points is defined to be the length
of the shortest connecting path:
$\Gdist{x}{y} = 2 \arcsin \tfrac{\Edist{x}{y}}{2}$.
This shortest path is unique, unless $y = -x$, in which case there are
infinitely many shortest paths of length $\pi$.
Letting $X$ be a finite set of points on $\Sspace^n$, we define
the \emph{Voronoi domain} of $x \in X$ as the points
for which $x$ minimizes the geodesic distance,
further constraining it to within the open hemisphere centered at $x$:
\begin{align}
  \Voronoi{x}  &=  \{ w \in \Sspace^n
    \mid  \Gdist{w}{x} \leq \Gdist{w}{y} \mbox{\rm ~for all~} y \in X 
          \mbox{\rm ~and~} \Gdist{w}{x} < \tfrac{\pi}{2} \}.
  \label{eqn:VorDef}
\end{align}
Note that $\Gdist{w}{x} \leq \Gdist{w}{y}$ defines a closed hemisphere,
namely all points $w \in \Sspace^n$ that satisfy
$\Edist{w}{x} \leq \Edist{w}{y}$ in $\Rspace^{n+1}$.
It follows that $\Voronoi{x}$ is the intersection of a finite collection
of hemispheres --- a set we refer to as a \emph{(convex) spherical polytope}.
Any two of these spherical polytopes have disjoint interiors.
The \emph{Voronoi tessellation} of $X$ is the collection of Voronoi domains,
one for each point in $X$.
It covers the entire $n$-sphere, except if $X$ is contained in a closed
hemisphere, in which case it covers $\Sspace^n$ minus a possibly
degenerate but non-empty spherical polytope;
see Figure \ref{fig:Voronoi}.
Generically, the common intersection of $1 \leq k \leq n+1$
Voronoi domains is either empty or a shared face of dimension $n-k+1$,
and the common intersection of $n+2$ or more Voronoi domains is empty.
The \emph{Delaunay mosaic} of $X$ is isomorphic to the nerve of the
Voronoi tessellation:
\begin{align}
  \Delaunay{}{X}  &=  \{ Q \subseteq X 
    \mid  \bigcap\nolimits_{x \in Q} \Voronoi{x} \neq \emptyset \} .
  \label{eqn:DelDef}
\end{align}
The Nerve Theorem \cite{Ler45} implies that the Delaunay mosaic has the same
homotopy type as the union of Voronoi domains.
Assuming there is no closed hemisphere that contains all points,
this is the homotopy type of $\Sspace^n$.

%%%%%%%%%%%%%%%%%%%%%%%%%%%%%%%%%%%%%%%%%%%%%%%%%%%%%%%%%%%%%%%%%%%%%%%%%%
\ourparagraph{Delaunay mosaics and inscribed polytopes.}
%%%%%%%%%%%%%%%%%%%%%%%%%%%%%%%%%%%%%%%%%%%%%%%%%%%%%%%%%%%%%%%%%%%%%%%%%%
The Delaunay mosaic is an (abstract) simplicial complex.
In the generic case, $\Delaunay{}{X}$ can be geometrically realized in
$\Rspace^{n+1}$, namely by mapping every abstract simplex, $Q$,
to the convex hull of its points.
To make this precise, we compare $\Delaunay{}{X}$ with the boundary complex
of $\conv{X}$, which is a convex polytope inscribed in the $n$-sphere.
Each $(n-1)$-sphere $S \subseteq \Sspace^n$ defines two (closed) caps.
If $S$ is a great-sphere, these caps are hemispheres,
else they have different volume and we call one the \emph{small cap}
and the other the \emph{big cap}.
Every facet of $\conv{X}$ defines such a pair of caps, namely the
portions of $\Sspace^n$ on the two sides of the $n$-plane spanned by the facet.
One of these caps is \emph{empty}, by which we mean that no point of $X$
lies in its interior.
If $0$ is in the interior of $\conv{X}$, then all empty caps are small,
but if $0 \not\in \conv{X}$, then there is at least one empty big cap.
For non-generic sets, $0$ may lie on the boundary of $\conv{X}$,
in which case there is at least one empty hemisphere cap.
Parsing the definitions in \eqref{eqn:VorDef} and \eqref{eqn:DelDef},
we observe that a simplex
$Q \subseteq X$ belongs to the Delaunay mosaic iff there is
an $(n-1)$-sphere, $S$, that contains $Q$, is not a great-sphere,
and whose empty cap is small.
In the generic case, these simplices $Q$ are exactly the faces of the
facets of $\conv{X}$ whose small caps are empty.
In particular, it shows that if points are not contained in any hemisphere,
then $\Delaunay{}{X}$ is isomorphic to $\conv{X}$, a random inscribed polytope.

%%%%%%%%%%%%%%%%%%%%%%%%%%%%%%%%%%%%%%%%%%%%%%%%%%%%%%%%%%%%%%%%%%%%%%%%%%
\ourparagraph{Radius function.}
%%%%%%%%%%%%%%%%%%%%%%%%%%%%%%%%%%%%%%%%%%%%%%%%%%%%%%%%%%%%%%%%%%%%%%%%%%
Consider growing a spherical cap from each point in $X$.
To formalize this process, we recall that
$\Capp{\GRad}{x} = \{ w \in \Sspace^n \mid \Gdist{w}{x} \leq \GRad \}$
is the cap with center $x \in X$ and geodesic radius $\GRad$.
Clipping the Voronoi domain to within the cap, for each point $x \in X$,
we get a subcomplex of the Delaunay mosaic when we take the nerve:
\begin{align}
  \Delaunay{\GRad}{X}  &=  \{ Q \subseteq X \mid
    \bigcap\nolimits_{x \in Q} [\Voronoi{x} \cap \Capp{\GRad}{x}] \neq \emptyset \}.
\end{align}
By construction, $\Delaunay{\GRad}{X}$ is a simplicial complex,
which we call the \emph{Delaunay complex},
and $\Delaunay{\GRad}{X} \subseteq \Delaunay{\GRadPr}{X}$
whenever $\GRad \leq \GRadPr$.
For $\GRad = \tfrac{\pi}{2}$, each restricting cap is a hemisphere
and thus contains its corresponding Voronoi domain, which implies
$\Delaunay{\pi/2}{X} = \Delaunay{}{X}$.
We are now ready to introduce the \emph{radius function},
$\Rfun \colon \Delaunay{}{X} \to \Rspace$,
which maps every simplex to the smallest geodesic radius for which
the simplex belongs to the subcomplex of the Delaunay mosaic:
\begin{align}
  \Rfun (Q)  &=  \min \{ \GRad \mid Q \in \Delaunay{\GRad}{X} \} .
\end{align}
In other words, $\Rfun^{-1} [0,\GRad] = \Delaunay{\GRad}{X}$.
This definition is different from but equivalent to
the one we gave in the introduction.
We will prove shortly that for generic $X$, the radius function on
the Delaunay mosaic is a generalized discrete Morse function;
see \cite{BaEd16,For98}.
To explain what this means, we let $L \subseteq U$ be two simplices
in $\Delaunay{}{X}$,
and we call $[L,U] = \{ Q \mid L \subseteq Q \subseteq U \}$ an \emph{interval}
and $(\ell, k)$ with $\ell = \dime{L}$ and $k = \dime{U}$ its \emph{type}.
For simple combinatorial reasons, the number of simplices in $[L,U]$
is $2^{k-\ell}$.
A function $g \colon \Delaunay{}{X} \to \Rspace$ is a
\emph{generalized discrete Morse function} if there exists a partition
of $\Delaunay{}{X}$ into intervals such that
$g(P) \leq g(Q)$ whenever $P \subseteq Q$,
with equality in this case iff $P$ and $Q$ belong to the same interval.
We can prove that the radius function for a generic
set $X$ satisfies this condition.
Formally, we say a finite set $X \subseteq \Sspace^n$
is in \emph{general position} if $\card{X} > n+1$ and for every $0 \leq k < n$
\medskip \begin{enumerate}
  \item[1.] no $k+3$ points of $X$ belong to a common $k$-sphere
            on $\Sspace^n$,
  \item[2.] considering the unique $(k+1)$-sphere that passes through
            $k+3$ points of $X$, no $k+2$ of these points belong to a
            common $k$-sphere that shares its center with the $(k+1)$-sphere.
\end{enumerate} \medskip
Condition 2 implies that no
$n+1$ points of $X$ lie on a great-sphere of $\Sspace^n$.
We need a few additional concepts.
Assume $X$ is in general position and $Q \subseteq X$ is a $k$-simplex
with $0 \leq k \leq n$.
A cap \emph{circumscribes} $Q$ if the bounding $(n-1)$-sphere passes
through all points of $Q$.
Since $X$ is generic, $Q$ has a unique \emph{smallest circumscribed cap},
which we denote $\capp{Q}$.
If $Q \in \Delaunay{}{X}$, $Q$ also has a unique \emph{smallest empty circumscribed cap},
which may or may not be the smallest circumscribed cap.
We call it the \emph{circumcap} of $Q$ and denote it as $\ecapp{Q}$.
The \emph{Euclidean center} of a cap is the center of the
bounding $(n-1)$ sphere, which
is a point in $\Rspace^{n+1}$ but not on $\Sspace^n$.
Using this center, we introduce a notion of visibility within the affine hull of $Q$,
which is a $k$-dimensional plane in $\Rspace^{n+1}$.
Recalling that a facet of $k$-simplex is a $(k-1)$-dimensional face,
we say a facet of $Q$ is \emph{visible} from this center if the
$(k-1)$-plane spanned by the facet separates the center from $Q$ or,
equivalently, if the center lies in one closed $k$-dimensional halfspace
bounded by the $(k-1)$-plane and $Q$ is contained in the other such halfspace.
\begin{lemma}[Radius Function]
  \label{lem:RadiusFunction}
  Let $X \subseteq \Sspace^n$ be finite and in general position.
  Then $\Rfun \colon \Delaunay{}{X} \to \Rspace$ is a generalized
  discrete Morse function,
  and $[L,U]$ is an interval of $\Rfun$ iff $\capp{U}$ is empty
  and $L$ is the maximal common face of all facets of $U$
  that are visible from the Euclidean center of $\capp{U}$.
  Furthermore, for every $Q \in [L, U]$, we have $\ecapp{Q} = \capp{U}$.
\end{lemma}
\ourproof
  We prove that for each $Q \in \Delaunay{}{X}$ there are unique
  Delaunay simplices $L \subseteq Q \subseteq U$ such that
  $\capp{U} = \ecapp{U}$,
  $L$ is the intersection of all visible facts of $U$,
  and all simplices in $[L,U]$ share the circumcap.
  Note that $\Rfun (Q)$ is the geodesic radius of the circumcap of $Q$.	
  Letting $U \subseteq X$ be the set of all points on the $(n-1)$-sphere
  that bounds this circumcap,
  we have $\ecapp{U} = \capp{U}$ for else we could find
  a smaller empty circumscribed cap.
  Let $z$ be the center and $\GRad$ the geodesic radius of $\capp{U}$.
  By assumption of general position,
  $|U| \leq n+1$, so $U$ is a Delaunay simplex.
  For every facet $F$ of $U$, let $z_F$ be the center and $\GRad_F$
  the geodesic radius of $\capp{F}$,
  and let $u_F$ be the unique vertex in $U \setminus F$.
  We move the center of this cap along the shortest path from $z_F$ to $z$
  while adjusting the radius so that all points of $F$ remain on the
  boundary of the cap.
  During this motion, the radius increases continuously, and when it
  reaches $\GRad$, the boundary of the cap passes through $u_F$.
  If $F$ is visible from $z$,
  then $u_F$ is inside the cap at the beginning and on the boundary
  of the cap at the end of the motion.
  If $F$ is not visible from the Euclidean center,
  then $u_F$ changes from outside at the beginning to on the boundary
  of the cap at the end of the motion.
  In other words, $\capp{U}$ is the circumcap of every
  visible facet of $U$,
  but every invisible facet has a smaller empty circumscribed cap.
  Since the intersection of two simplices with common circumcap
  has the same circumcap
  \cite[Lemma 3.4]{BaEd16},
  we can take $L$ as the intersection of all visible facets of $U$
  and get $\ecapp{L} = \capp{U}$.
  On the other hand, any face of $U$ that does not contain $L$
  is also a face of an invisible facet
  and therefore has a smaller empty circumscribed cap.
  This implies $L \subseteq Q$.
	
  We note that the construction gives a partition of $\Delaunay{}{X}$
  into intervals.
  Indeed, any two Delaunay simplices sharing the circumcap
  give rise to the same simplex $U$ and therefore to the same
  interval $[L, U]$.
  This concludes the proof.
\eop

\noindent {\sc Remark.}
({\sc 4a}) While the proof follows almost verbatim the proof
in the Euclidean case \cite{BaEd16},
and actually the Euclidean Delaunay mosaic of the spherical point set
is almost identical to the one we are talking about,
there is a subtlety hidden in its definition.
Indeed, because each Voronoi domain is restricted to within the
open hemisphere centered at the generating point, the sets
$\Voronoi{x} \cap \Capp{\GRad}{x}$ form a system
in which every common intersection is either empty or contractible.
The Nerve Theorem thus applies, proving that the
subcomplex of the Delaunay mosaic has the same homotopy type as the union
of caps of radius $\GRad$.
This property breaks down for the boundary complex of $\conv{X}$.
This can be seen by considering the four points on $\Sspace^2$ shown
in Figure \ref{fig:Voronoi}:
$A, B = (\pm \ee, 0, \sqrt{1-\ee^2})$ and $C, D = (0, \pm 1/2, \sqrt{3}/2)$,
in which $\ee$ is a sufficiently small positive real number.
The great-circle arc shared by the Voronoi domains of $C$ and $D$
has length only slightly shorter than $\pi$
and it intersects the union of four caps
of geodesic radius $\GRad$ slightly larger than $\tfrac{\pi}{2}$
in two disconnected segments.
The union of the four caps has the topology of a disk,
while the nerve has the topology of a circle.
Indeed, the latter consists of two triangles glued along a shared edge
plus another edge connecting the two respective third vertices of the two triangles.

%%%%%%%%%%%%%%%%%%%%%%%%%%%%%%%%%%%%%%%%%%%%%%%%%%%%%%%%%%%%%%%%%%%%%%%%%%
\ourparagraph{Poisson point process.}
%%%%%%%%%%%%%%%%%%%%%%%%%%%%%%%%%%%%%%%%%%%%%%%%%%%%%%%%%%%%%%%%%%%%%%%%%%
We are interested in sets $X \subseteq \Sspace^n$ that are randomly generated.
In particular, we use a (stationary) \emph{Poisson point process} with \emph{density}
$\density > 0$, which is characterized by the following two properties:
\medskip \begin{enumerate}
  \item[1.]  the numbers of points in a finite collection of
    pairwise disjoint Borel sets on $\Sspace^n$ are independent random
    variables;
  \item[2.]  the expected number of points in a Borel set is
    $\density$ times the Lebesgue measure of the set.
\end{enumerate} \medskip
See \cite{Kin93} for an introduction to Poisson point processes.
The two conditions imply that the probability of having $k$ points in
a Borel set $B \subseteq \Sspace^n$ with Lebesgue measure $\norm{B}$ is
$\Probable{\card{X \cap B} = k}
   = \density^k \norm{B}^k e^{- \density \norm{B}} / k!$.
In particular, the probability of having no point in $B$ is
$\Probable{X \cap B = \emptyset} = e^{- \density \norm{B}}$.
It is not difficult to prove that the realization $X$ of
a Poisson point process on $\Sspace^n$ is finite and in general position
with probability $1$, a property we will assume for the remainder
of this paper.
It follows that $\Delaunay{}{X}$ is an $n$-dimensional simplicial complex
and, by Lemma \ref{lem:RadiusFunction}, that
$\Rfun \colon \Delaunay{}{X} \to \Rspace$ is a generalized
discrete Morse function.

To familiarize ourselves with the definition of a Poisson point process, we prove that
the difference between the boundary complex of $\conv{X}$ and
$\Delaunay{}{X}$ is small.
More precisely, the number of faces of $\conv{X}$ that are visible
from $0$ outside $\conv{X}$ vanishes rapidly as the density increases.
This is consistent with the rapid decrease of the probability that
$0 \not\in \conv{X}$, as computed by Wendel \cite{Wen62}
for the uniform distribution on $\Sspace^n$.
\begin{lemma}[Non-Delaunay Faces]
  \label{lem:NonDelaunayFaces}
  Let $X$ be a Poisson point process with density $\density > 0$
  on $\Sspace^n$.
  For every $0 \leq k \leq n$, the expected number of $k$-faces
  of $\conv{X}$ that do not belong to $\Delaunay{}{X}$ goes to $0$
  as $\density$ goes to $\infty$.
\end{lemma}
\ourproof
  We may assume that $\conv{X}$ is simplicial and that no $n+1$ points
  lie on a great-sphere of $\Sspace^n$.
  Let $Q \subseteq X$ be a set of $n+1$ points and consider its
  small and big caps.
  The big cap has volume larger than $\sigma_{n+1} / 2$,
  and $Q$ is a facet of $\conv{X}$ but not a simplex of $\Delaunay{}{X}$
  iff this big cap is empty.
  The probability of this event is less than $e^{- \density \sigma_{n+1}/2}$.
  The expected number of such facets of $\conv{X}$ is therefore less than
  a constant times $\density^{n+1} e^{- \density \sigma_{n+1} / 2}$,
  which goes to $0$ as $\density$ goes to $\infty$.
  Here we used that $\Expected{|X|^{n+1}}$
  is at most a constant times $\density^{n+1}$.
  For $k < n$, every $k$-face of $\conv{X}$ that does not belong
  to $\Delaunay{}{X}$ is a face of a facet with this property.
  The expected number of such $k$-faces thus also goes to $0$
  as $\density$ goes to $\infty$.
\eop

%\newpage
%%%%%%%%%%%%%%%%%%%%%%%%%%%%%%%%%%%%%%%%%%%%%%%%%%%%%%%%%%%%%%%%%%%%%%%%%%
%%%%%%%%%%%%%%%%%%%%%%%%%%%%%%%%%%%%%%%%%%%%%%%%%%%%%%%%%%%%%%%%%%%%%%%%%%
\section{Proof of Main Result}
\label{sec:4}
%%%%%%%%%%%%%%%%%%%%%%%%%%%%%%%%%%%%%%%%%%%%%%%%%%%%%%%%%%%%%%%%%%%%%%%%%%
%%%%%%%%%%%%%%%%%%%%%%%%%%%%%%%%%%%%%%%%%%%%%%%%%%%%%%%%%%%%%%%%%%%%%%%%%%

In this section, we prove the main result of this paper stated
as Theorem \ref{thm:MainResult} in the Introduction.
It consists of an integral equation for the expected number of intervals
as a function of the maximum geodesic radius,
and an asymptotic version of the formula for $\density \to \infty$.

%%%%%%%%%%%%%%%%%%%%%%%%%%%%%%%%%%%%%%%%%%%%%%%%%%%%%%%%%%%%%%%%%%%%%%%%%%
%%%%%%%%%%%%%%%%%%%%%%%%%%%%%%%%%%%%%%%%%%%%%%%%%%%%%%%%%%%%%%%%%%%%%%%%%%
\subsection{The Integral Equation}
\label{sec:41}
%%%%%%%%%%%%%%%%%%%%%%%%%%%%%%%%%%%%%%%%%%%%%%%%%%%%%%%%%%%%%%%%%%%%%%%%%%
%%%%%%%%%%%%%%%%%%%%%%%%%%%%%%%%%%%%%%%%%%%%%%%%%%%%%%%%%%%%%%%%%%%%%%%%%%

We begin with the proof of the integral equation, \eqref{eqn:MainResult}.
The main tools are the Slivnyak--Mecke formula,
which we will discuss shortly,
and the Blaschke--Petkantschin formula for the sphere,
which was stated and proved in Section \ref{sec:2}.
In addition, we employ the combinatorial analysis
of inscribed simplices in \cite{ENR16}.

%%%%%%%%%%%%%%%%%%%%%%%%%%%%%%%%%%%%%%%%%%%%%%%%%%%%%%%%%%%%%%%%%%%%%%%%
\ourparagraph{The Slivnyak--Mecke approach.}
%%%%%%%%%%%%%%%%%%%%%%%%%%%%%%%%%%%%%%%%%%%%%%%%%%%%%%%%%%%%%%%%%%%%%%%%
In a nutshell, the \emph{Slivnyak--Mecke formula} writes the expectation
of a random variable of a Poisson point process as an integral over
the space on which the process is defined; see \cite[page 68]{ScWe08}.
To write this integral, we recall that $\xxx = (x_0, x_1, \ldots, x_k)$
is a sequence of $k+1$ points or $k$-simplex on $\Sspace^n$,
that $\PemptyOnly \colon \left( \Sspace^n \right)^{k+1} \to \Rspace$
maps $\xxx$ to the probability that its smallest circumscribed cap is empty,
that $\One{k-\ell} \colon (\Sspace^n)^{k+1} \to \Rspace$
indicates whether or not the number of facets visible from
the Euclidean center of the smallest circumscribed cap is $k - \ell$,
and that $\One{\GRad} \colon \left( \Sspace^n \right)^{k+1} \to \Rspace$
indicates whether or not $\Rfun (\xxx) \leq \GRad$.
Choosing points from a Poisson point process with density $\density > 0$
on $\Sspace^n$, we use Slivnyak--Mecke to write the expected number
of intervals of type $(\ell, k)$ and geodesic radius at most $\GRad_0$ as
\begin{align}
  \Expected{\ccon{\ell}{k}{n}, \GRad_0}
    &=  \tfrac{\density^{k+1}}{(k+1)!}
        \int\displaylimits_{\xxx \in (\Sspace^n)^{k+1}}
        \Pempty{\xxx} \cdot \One{k-\ell} (\xxx) \cdot
                            \One{\GRad_0} (\xxx) \diff \xxx ,
  \label{eqn:SM}
\end{align}
in which $0 \leq \ell \leq k \leq n$; compare with \cite{ENR16}.
The probability that the smallest circumscribed cap of the $k$-simplex is empty
is $\Pempty{\xxx} = e^{- \density \Area{\GRad}}$,
with $\GRad$ the geodesic radius of the cap.
To compute the integral in \eqref{eqn:SM}, we apply
Equation \eqref{eqn:rotationally-symmetric} from Theorem \ref{thm:BPforSphere}
with $f(\xxx) = \Pempty{\xxx} \One{k-\ell} (\xxx)\One{\GRad_0} (\xxx)$.
The corresponding function from the statement of Theorem \ref{thm:BPforSphere},
$f_r\colon (\Sspace^{k-1})^{k+1} \subseteq (\Rspace^{n+1})^{k+1} \to \Rspace$,
is defined by
$f_r (\uuu) = \Pempty{r} \One{k-\ell} (\uuu) \One{\GRad_0} (r)$,
where we write $\Pempty{r} = \Pempty{\uuu}$ and $\One{\GRad_0}(r) = \One{\GRad_0}(\uuu)$
to emphasize that these expressions depend only on the radius.
Equation \eqref{eqn:rotationally-symmetric} then gives
\begin{align}
  \int\displaylimits_{\xxx \in (\Sspace^n)^{k+1}} f(\xxx) \diff \xxx
    &=  \tfrac{\sigma_{n+1}}{2} \norm{\LGrass{k}{n}} 
        \int\displaylimits_{t = 0}^1 t^{\frac{kn-2}{2}} (1-t)^{\frac{n-k-1}{2}}
          \Pempty{\sqrt{t}} \One{\GRad_0} (\sqrt{t})           \nonumber \\
      & \quad\quad\quad
        \times \int\displaylimits_{\uuu \in (\Sspace^{k-1})^{k+1}}
          \One{k-\ell} (\uuu) [k!\Volume{\uuu}]^{n-k+1} \diff \uuu \diff t.
    \label{eqn:tryit}
\end{align}

%%%%%%%%%%%%%%%%%%%%%%%%%%%%%%%%%%%%%%%%%%%%%%%%%%%%%%%%%%%%%%%%%%%%%%%%
\ourparagraph{Substitution and reformulation.}
%%%%%%%%%%%%%%%%%%%%%%%%%%%%%%%%%%%%%%%%%%%%%%%%%%%%%%%%%%%%%%%%%%%%%%%%
To continue, we recall
$\EAux{\ell}{k}{n} = \Expected{\Volume{\uuu}^{n-k+1} \One{k-\ell}(\uuu)}$
from \eqref{eqn:constantE},
in which the expectation is for sampling $k+1$ points from
the uniform distribution on $\Sspace^{k-1}$.
It follows that the second integral on the right-hand side of
\eqref{eqn:tryit} is $k!^{n-k+1} \sigma_k^{k+1} \EAux{\ell}{k}{n}$.
Rewriting \eqref{eqn:SM} using \eqref{eqn:tryit}, we therefore get
\begin{align}
  \Expected{\ccon{\ell}{k}{n},\GRad_0}  &= 
    \tfrac{\density^{k+1}}{(k+1)!}
      \tfrac{\sigma_{n+1}}{2} \norm{\LGrass{k}{n}}
      k!^{n-k+1} \sigma_k^{k+1} \EAux{\ell}{k}{n}
      \int\displaylimits_{t = 0}^{s}
        t^{\frac{kn-2}{2}} (1-t)^{\frac{n-k-1}{2}} \Pempty{\sqrt{t}} \diff t,
  \label{eqn:f10first} \\
    &=  \density \sigma_{n+1}
      \cdot \tfrac{\sigma_n^k}{2 \Gamma(k) n^{k-1}}
      \cdot \Ccon{\ell}{k}{n}
      \int\displaylimits_{t = 0}^{s}
        \density^k t^{\frac{kn-2}{2}} (1-t)^{\frac{n-k-1}{2}} \Pempty{\sqrt{t}} \diff t,
  \label{eqn:f10last}
\end{align}
in which we absorb one indicator by limiting the range of integration
to the square of the maximum Euclidean radius,
$s = \sin^2 \GRad_0$.
To get \eqref{eqn:f10last} from \eqref{eqn:f10first}, we cancel $k!$,
move $\density^k$ inside the integral,
and use \eqref{eqn:constantE} to substitute
$[\sigma_n^k / (\Gamma(k) n^{k-1})] \cdot \Ccon{\ell}{k}{n}$ for
$[\norm{\LGrass{k}{n}} k!^{n-k} \sigma_k^{k+1} / (k+1)] \cdot \EAux{\ell}{k}{n}$.
This proves the integral equation \eqref{eqn:MainResult} in Theorem \ref{thm:MainResult}.
%% We will see shortly that after moving $\density^k$, the integral is bounded
%% away from $0$ as well as from $\infty$.
%% It follows that the expected number of intervals grows linearly with the density.

%%%%%%%%%%%%%%%%%%%%%%%%%%%%%%%%%%%%%%%%%%%%%%%%%%%%%%%%%%%%%%%%%%%%%%%%%%
%%%%%%%%%%%%%%%%%%%%%%%%%%%%%%%%%%%%%%%%%%%%%%%%%%%%%%%%%%%%%%%%%%%%%%%%
\subsection{The Asymptotic Result}
\label{sec:42}
%%%%%%%%%%%%%%%%%%%%%%%%%%%%%%%%%%%%%%%%%%%%%%%%%%%%%%%%%%%%%%%%%%%%%%%%%%
%%%%%%%%%%%%%%%%%%%%%%%%%%%%%%%%%%%%%%%%%%%%%%%%%%%%%%%%%%%%%%%%%%%%%%%%%%
We continue with the proof of the asymptotic result \eqref{eqn:asymp_res}.
We proceed in two stages,
first taking liberties and leaving gaps in the argument,
and second filling all the gaps.

%%%%%%%%%%%%%%%%%%%%%%%%%%%%%%%%%%%%%%%%%%%%%%%%%%%%%%%%%%%%%%%%%%%%%%%%%%%%
\ourparagraph{Argument with gaps.}
%%%%%%%%%%%%%%%%%%%%%%%%%%%%%%%%%%%%%%%%%%%%%%%%%%%%%%%%%%%%%%%%%%%%%%%%%%%%
We are interested in the behavior of the integral in \eqref{eqn:f10last},
when $\density \to \infty$.
We observe that the probability of a cap to be empty
vanishes rapidly with increasing geodesic radius:
$\Pempty{\ERad} = e^{- \density \Area{\GRad}}$,
in which $\ERad = \sin \GRad$ is the Euclidean radius.
This implies that the integrand is concentrated in the vicinity of $0$.
To make sense of the radius in the limit,
we re-scale by mapping $\GRad$ and $\density$
to the normalized radius, $\NRad = \GRad \density^{1/n}$.
To proceed with the informal computations, we assume that $\GRad$
is close to $0$ and prepare two approximations and one relation:
\medskip \begin{description}
  \item[{\rm A.}]  the squared Euclidean radius is roughly the squared
    geodesic radius:  $s = \sin^2 \GRad \approx \GRad^2$;
  \item[{\rm B.}]  the square of the height is $1 - s \approx 1$,
    which allows us to simplify the incomplete Beta function:
    \begin{align}
      \iBeta{s}{\tfrac{n}{2}}{\tfrac{1}{2}}
        &=       \int_{t=0}^s
                   t^{\frac{n}{2} - 1} (1-t)^{- \frac{1}{2}} \diff t
        \approx  \int_{t=0}^s t^{\frac{n}{2} - 1} \diff t
        =        \tfrac{2}{n} s^{n/2} ;
    \end{align}
  \item[{\rm C.}]  the relation $\tfrac{\sigma_{n+1}}{\sigma_n} = \Beta{\tfrac{n}{2}}{\tfrac{1}{2}}$
    implies $\tfrac{\sigma_{n+1}}{n} / \Beta{\tfrac{n}{2}}{\tfrac{1}{2}}
              = \tfrac{\sigma_n}{n} = \nu_n$.
\end{description} \medskip
Returning to the integral in \eqref{eqn:MainResult},
but without the factor $\density^n$, we get
\begin{align}
  \int_{t=0}^{\sin^2 \GRad_0}
    t^{\frac{kn-2}{2}} (1-t)^{\frac{n-k-1}{2}}
    \Pempty{\sqrt{t}} \diff t
  &\approx \int_{t=0}^{\NRad_0^2/\density^{2/n}}
    t^{\frac{kn-2}{2}} e^{- \density \nu_n t^{n/2}} \diff t ,
  \label{eqn:intuit1}
\end{align}
in which we approximate the upper limit of the integration using A,
and drop the middle factor because it is close to $1$ according to B.
The probability of having an empty cap is
$\Pempty{\ERad} = e^{- \density \Area{\GRad}}$,
in which the area of the cap can be written in terms of Beta functions:
\begin{align}
  \Area{\GRad}  &=  \tfrac{\sigma_{n+1} \iBeta{s}{n/2}{1/2}}
                          {2 \Beta{n/2}{1/2}}
    \approx  \tfrac{\sigma_{n+1} (2/n) s^{n/2}}
                   {2 \Beta{n/2}{1/2}}
    =        \nu_n s^{n/2} ,
  \label{eqn:intuit1b}
\end{align}
using B for the approximation and C to get the final result,
which we plug into the left-hand side of \eqref{eqn:intuit1}
to get the approximation on its right-hand side.
The exponential term motivates us to change variables
with $\tau = \density \nu_n t^{n/2}$.
Plugging $t = \tau^{2/n} / (\density \nu_n)^{2/n}$
and $\diff t = [ \tfrac{2}{n} \tau^{2/n-1} / (\density \nu_n)^{2/n}] \diff \tau$
into the right-hand side of \eqref{eqn:intuit1}, we get
\begin{align}
  \int_{\tau=0}^{v} \tau^{k-1} (\density \nu_n)^{-k}
                  \left( \tfrac{2}{n} \right) e^{- \tau} \diff \tau
    &=  \tfrac{2 n^{k-1}}{\density^k \sigma_n^k} \cdot \iGama{v}{k} ,
  \label{eqn:intuit2}
\end{align}
in which the upper bound of the integration range is
$v = \density \nu_n ( \NRad_0^2 / \density^{2/n} )^{n/2} = \NRad_0^n \nu_n$,
the power of $\tau$ is $\tfrac{2}{n} \tfrac{kn-2}{2} + \tfrac{2}{n} - 1 = k-1$,
and the power of $\density \sigma_n$ is
$- \tfrac{2}{n} \tfrac{kn-2}{2} - \tfrac{2}{n} = - k$.
We get the right-hand side of \eqref{eqn:intuit2} from the
left-hand side using $\tfrac{\sigma_n}{n} = \nu_n$
and $\iGama{v}{k} = \int_{\tau=0}^v \tau^{k-1} e^{- \tau} \diff \tau$.
Finally plugging the right-hand side into \eqref{eqn:MainResult}, we get
\begin{align}
  \Expected{\ccon{\ell}{k}{n}, \GRad_0} &= 
    \density \sigma_{n+1} \cdot \tfrac{\sigma_n^k}{2 \Gamma(k) n^{k-1}}
      \cdot \Ccon{\ell}{k}{n}
      \int_{t = 0}^{\sin^2 \GRad_0}
        \density^k t^{\frac{kn-2}{2}} (1-t)^{\frac{n-k-1}{2}}
        \Pempty{\sqrt{t}} \diff t                              \\
    &=  \density \sigma_{n+1} \cdot 
        \tfrac{\iGama{v}{k}}{\Gama{k}} \cdot
        \Ccon{\ell}{k}{n} + \littleoh{\density} ,
  \label{eqn:intuit3}
\end{align}
as claimed in Theorem \ref{thm:MainResult}.
Making the unjustified substitution $v = \NRad_0^n \nu_n = \infty$, we get
\begin{align}
  \Expected{\ccon{\ell}{k}{n}}
    &=  \density \sigma_{n+1} \cdot  \Ccon{\ell}{k}{n} + \littleoh{\density} ,
  \label{eqn:intuit4}
\end{align}
as claimed in Remark ({\sc 1c}) after Theorem \ref{thm:MainResult}.

%%%%%%%%%%%%%%%%%%%%%%%%%%%%%%%%%%%%%%%%%%%%%%%%%%%%%%%%%%%%%%%%%%%%%%%%%%%%
\ourparagraph{Formal justifications.}
%%%%%%%%%%%%%%%%%%%%%%%%%%%%%%%%%%%%%%%%%%%%%%%%%%%%%%%%%%%%%%%%%%%%%%%%%%%%
We continue with the justification of the asymptotic equivalences
claimed above.
To recall, there is the approximation in \eqref{eqn:intuit1}
and the substitution $\NRad_0 = \infty$ after \eqref{eqn:intuit3}.
Fixing a real number $0 \leq \dd \leq 1$,
we introduce some notation to streamline the computations:
\begin{align}
  \alpha  &=  \tfrac{kn-2}{2},  \quad \alpha' = \tfrac{n-k-1}{2}, \quad
    \beta  =  \tfrac{n}{2},     \quad \beta'  = \tfrac{1}{2},     \quad
    c      =  \tfrac{\sigma_n}{2},                                         \\
  g(s)    &=  c \int_{t=0}^{s}
              t^{\beta-1} (1-t)^{\beta'-1} \diff t ,     \label{eqn:g}     \\
  J_0     &=  \density^k \int_{t=0}^1
              t^\alpha (1-t)^{\alpha'}e^{-\density g(t)} \diff t,  \quad
  J_1(\dd) =  \density^k \int_{t=0}^\dd
              t^\alpha (1-t)^{\alpha'}e^{-\density g(t)} \diff t ,         \\
  J_2(\dd) &= \density^k \int_{t=0}^\dd
              t^\alpha e^{-\density g(t)} \diff t , \quad \quad \quad \quad ~~
  J_3(\dd) =  \density^k \int_{t=0}^\dd
              t^\alpha e^{-\density \frac{c}{\beta}t^\beta} \diff t.
\end{align}
We note that $\alpha, \alpha' \geq - \tfrac{1}{2}$,
$\beta, \beta' \geq \tfrac{1}{2}$,
and $g(s)$ is $c = \tfrac{\sigma_n}{2}$ times the incomplete Beta function.
Recall that
$\tfrac{\sigma_{n+1}}{\sigma_n} = \Beta{\tfrac{n}{2}}{\tfrac{1}{2}}$,
which implies that $g(s)$ is $\tfrac{\sigma_{n+1}}{2}$ times the ratio of
incomplete over complete Beta functions.
Hence $g(s) = \Area{\GRad}$, in which $s = \sin^2 \GRad$;
see \eqref{eqn:CapArea}.
Note also that $J_0$ is the integral in \eqref{eqn:f10last} except that the
integration range goes all the way to $1$,
which corresponds to computing the number of intervals without
restricting the radius.
For $\dd = 1$, we have $J_1 = J_0$,
and for $\dd = \sin^2 \GRad_0$, $J_1$ is $\density^k$ times
the expression on the left-hand side of \eqref{eqn:intuit1}.
Finally, for $\dd = \NRad_0^2 \density^{-1/\beta}$,
$J_3$ is the integral on the right-hand side of \eqref{eqn:intuit1},
which we computed in \eqref{eqn:intuit2}.
Next, we list a sequence of observations:
\medskip \begin{description}
  \item[{\rm I.}]   The integral in \eqref{eqn:g} satisfies
    \begin{align}
      \tfrac{c}{\beta} s^\beta  &\leq  g(s)
        =     c \int_{t=0}^s t^{(n-2)/2} (1-t)^{-1/2} \diff t
        \leq  \tfrac{c}{\beta} s^\beta + \const \cdot s^{\beta+1},
    \end{align}
    for $0 \leq s \leq 1$ on the left,
    and for $0 \leq s \leq \tfrac{1}{2}$ on the right.
    Indeed, we have $1 \leq 1 / \sqrt{1-t}$ for all $0 \leq t \leq 1$ and
    $1 / \sqrt{1-t} \leq 1 + \const \cdot t$ for all $0 \leq t \leq \tfrac{1}{2}$.

  \medskip \item[{\rm II.}]
    The absolute difference between $J_0$ and $J_1(\dd)$ satisfies
    \begin{align}
      |J_0 - J_1(\dd)|  &=  \density^k \int_{t=\dd}^1 t^\alpha (1-t)^{\alpha'}
                                           e^{-\density g(t)} \diff t
                      \leq  \density^k e^{-\density \frac{c}{\beta} \dd^\beta}
                            \Beta{\alpha+1}{\alpha'+1} ,
    \end{align}
    because $g(t) \geq g(\dd)$ throughout the integration domain,
    and $g(\dd) \geq \tfrac{c}{\beta} t^\beta$ by I.
    The value of the Beta function is a constant independent of $\density$.

  \medskip \item[{\rm III.}]
    For $\dd \leq \tfrac{1}{2}$, the absolute difference
    between $J_1$ and $J_2$ satisfies
    \begin{align} 
      |J_1(\dd) - J_2(\dd)| 
        &\leq  \density^k \int_{t=0}^\dd [t^\alpha (1-t)^{\alpha'} - t^\alpha]
                 e^{- \density g(t)} \diff t
         \leq \const \cdot \delta J_2(\delta) ,
    \end{align}
    because $| 1 - (1-t)^{\alpha'} | \leq \const \cdot t$ for all
    $0 \leq t \leq \tfrac{1}{2}$ and $\alpha' \geq - \tfrac{1}{2}$.

  \medskip \item[{\rm IV.}]
    For $\dd \leq \tfrac{1}{2}$, the absolute difference
    between $J_2$ and $J_3$ satisfies
    \begin{align}
      |J_2(\dd) - J_3(\dd)|  &=  \density^k \int_{t=0}^\dd t^\alpha
        \left[ e^{-\density \frac{c}{\beta} t^\beta}
             - e^{-\density g(t)} \right] \diff t
          \label{eqn:IV1} \\
                          &\leq  \density^k \int_{t=0}^\dd t^\alpha
         e^{-\density \frac{c}{\beta} t^\beta}
        \left[ 1 - e^{-\const \cdot \density t^{\beta+1}} \right] \diff t
          \label{eqn:IV4} \\
                          &\leq J_3(\dd)
        \left[ 1 - e^{-\const \cdot \density \delta^{\beta+1}} \right],
          \label{eqn:IV5}
    \end{align}
    in which we use the left inequality in I to get the right order
    of the exponential terms in \eqref{eqn:IV1},
    and the right inequality in I to get \eqref{eqn:IV4}.

  \medskip \item[{\rm V.}]
    For $\GRad \leq 1 / \sqrt{2}$, the absolute difference between
    $J_1$ at the values $\sin^2 \GRad$ and $\GRad^2$ satisfies
    \begin{align}
      |J_1(\sin^2 \GRad) - J_1(\GRad^2)|
        &= \density^k \int_{t=\sin^2 \GRad}^{\GRad^2}
           t^\alpha (1-t)^{\alpha'} e^{-\density g(t)} \diff t
         \leq  2 \density^k \int_{t=\sin^2 \GRad}^{\GRad^2}
           t^\alpha \diff t
          \label{eqn:V2}    \\
	&\leq  \tfrac{2 \density^k}{\alpha+1}
            [ \GRad^{2\alpha + 2} - (\GRad - \GRad^2)^{2 \alpha +2} ]
         \leq  4 \density^k \GRad^{2\alpha + 3} ,
          \label{eqn:V5}
    \end{align}
    in which we use $(1-t)^{\alpha'} \leq 2$ for $t \leq \tfrac{1}{2}$
    to get the right-hand side of \eqref{eqn:V2}.
    We use $\sin \GRad > \GRad - \GRad^2$,
    which we glean from the Taylor series
    $\sin \GRad = \GRad - \tfrac{1}{6} \GRad^3 + \ldots$,
    and the binomial expansion of $(\GRad - \GRad^2)^{2 \alpha +2}$
    to get \eqref{eqn:V5}.
\end{description} \medskip

\noindent
As mentioned earlier, $J_1 (\sin^2 \GRad_0)$ is $\density^k$ times the
left-hand side of \eqref{eqn:intuit1},
and $J_3(\GRad_0^2)$ is $\density^k$ times the right-hand side
of \eqref{eqn:intuit1}.
According to \eqref{eqn:intuit2}, $\density^k$ times this right-hand side
is $(2 n^{k-1} / \sigma_n^k) \cdot \iGama{v}{k}$,
with $v = \NRad^n \nu_n$,
which is a positive constant;
see Remark ({\sc 1b}) where we first mentioned that this integral
is bounded from $0$ as well as from $\infty$.
Having established that there is a positive constant $C = J_3 (\GRad_0^2)$,
IV implies that
$J_2(\GRad_0^2) \leq C + (1-e^{- \density \frac{c}{\beta} \GRad_0^{2 (\beta+1)}})C$
is also bounded by a constant.
Using III, IV, V, we get
\begin{align}
  |J_1(\sin^2 \GRad_0) - J_3(\GRad_0^2)|
    &\leq  |J_1(\sin^2 \GRad_0) \!-\! J_1(\GRad_0^2)|
         + |J_1(\GRad_0^2) \!-\! J_2(\GRad_0^2)|
         + |J_2(\GRad_0^2) \!-\! J_3(\GRad_0^2)|               \\
    &\leq  4 \density^k \GRad_0^{2 \alpha + 3}
         + \const \cdot \GRad_0^2 J_2(\GRad_0^2)
         + (1 - e^{- \const \cdot \density \GRad_0^{2 (\beta+1)}} ) C .
    \label{eqn:ItoV}
\end{align}
Letting $\density$ to to infinity, we observe
\begin{align}
  \density^k \GRad_0^{2 \alpha + 3} &= \density^k \left(\NRad_0 \density^{-\tfrac{1}{n}}\right)^{kn+1} \to 0, \\
	\density \GRad_0^{2(\beta+1)} &= \density \left(\NRad_0 \density^{-\tfrac{1}{n}}\right)^{n+2} \to 0,
\end{align}
implying the three terms in \eqref{eqn:ItoV} go to $0$.
This finally justifies the approximation \eqref{eqn:intuit1} and the
argument proving Theorem \ref{thm:MainResult}.

%%%%%%%%%%%%%%%%%%%%%%%%%%%%%%%%%%%%%%%%%%%%%%%%%%%%%%%%%%%%%%%%%%%%%%%%%%%%
\ourparagraph{Justification of Remark ({\sc 1c}).}
%%%%%%%%%%%%%%%%%%%%%%%%%%%%%%%%%%%%%%%%%%%%%%%%%%%%%%%%%%%%%%%%%%%%%%%%%%%%
We finally prove that we can compute $J_0$ by setting
$\NRad_0$ to infinity in \eqref{eqn:intuit3} or, more formally,
by replacing the incomplete gamma function in the expression
by the complete gamma function.
Such a justification is needed because so far we have treated the
geodesic radius as a constant in our computations.
We now couple the bound of the integration domain with the density
by setting $\dd_0 = \density^{-1/(\beta+1/2)}$.
We reuse Equations \eqref{eqn:intuit1} and \eqref{eqn:intuit3} to
compute $J_3(\dd_0) = (2 n^{k-1} / \sigma_n^k) \cdot \iGama{v}{k}$,
with $v = \density \nu_n \dd_0^{n/2} = \nu_n \density^{1/(n+1)}$.
The upper bound for the incomplete Gamma function thus goes to infinity
and approaches the complete Gamma function.
We still have $J_3(\dd_0)$ bounded by a constant,
so the rest of the argument above goes through. We finally use II, which shows $|J_0- J_1(\dd_0)| \to 0$.
This justifies \eqref{eqn:intuit4} and Remark ({\sc 1c})
in the Introduction.

%\newpage
%%%%%%%%%%%%%%%%%%%%%%%%%%%%%%%%%%%%%%%%%%%%%%%%%%%%%%%%%%%%%%%%%%%%%%%%%%
%%%%%%%%%%%%%%%%%%%%%%%%%%%%%%%%%%%%%%%%%%%%%%%%%%%%%%%%%%%%%%%%%%%%%%%%%%
\section{Discussion}
\label{sec:5}
%%%%%%%%%%%%%%%%%%%%%%%%%%%%%%%%%%%%%%%%%%%%%%%%%%%%%%%%%%%%%%%%%%%%%%%%%%
%%%%%%%%%%%%%%%%%%%%%%%%%%%%%%%%%%%%%%%%%%%%%%%%%%%%%%%%%%%%%%%%%%%%%%%%%%

The main result of this paper is a radius-dependent integral equation for
the expected number of intervals of the radius function of
a Poisson point process on $\Sspace^n$.
To first order, the expected numbers are the same as in $\Rspace^n$;
compare with \cite{ENR16}.
The Delaunay mosaics on $\Sspace^n$ relate to inscribed convex polytopes
in $\Rspace^{n+1}$ and to the Delaunay mosaics in the standard $n$-simplex
equipped with the Fisher information metric.
These diagrams have therefore very similar stochastic properties
as the Delaunay mosaics in $\Rspace^n$.
We formulate a few questions that are motivated by the findings reported
in this article.
\medskip \begin{itemize}
  \item As mentioned earlier, the first-order terms of the expected
    number of intervals of the radius function do not
    distinguish $\Sspace^n$ from $\Rspace^n$.
    There are no further terms in the Euclidean case,
    but what are they for $\Sspace^n$?
  \item Projecting the convex hull of a finite $X \subseteq \Sspace^n$
    orthogonally onto a $(k+1)$-plane corresponds to slicing the Voronoi
    tessellation of $X$ with a $k$-dimensional great-sphere of $\Sspace^n$.
    Similarly, we can define a $k$-dimensional weighted Delaunay mosaic
    by slicing a Voronoi tessellation in $\Rspace^n$ with a $k$-plane.
    What are the stochastic properties of these slices and projections?
  \item The square of the Fisher information metric agrees infinitesimally
    with the Kullback--Leibler divergence \cite{KuLe51}.
    The more general class of Bregman divergences has recently come
    into focus \cite{EdWa16}.
    What are the stochastic properties of the Bregman divergences and their
    corresponding metrics?
    Is the similarity to the Euclidean metric specific to the
    Fisher information metric or is it a more general phenomenon?
\end{itemize}

\subsection*{Acknowledgements}
{\small The authors thank Matthias Reitzner for sharing a draft
  on Poisson--Delaunay mosaics on the sphere,
  {\v Z}iga Virk and Hubert Wagner for their help in connecting
  this work with Fisher information space,
  and Nicholas Barton for pointing out that the connection has been
  discovered earlier by Antonelli.}

%\newpage
%%%%%%%%%%%%%%%%%%%%%%%%%%%

\newpage \appendix
%%%%%%%%%%%%%%%%%%%%%%%%%%%%%%%%%%%%%%%%%%%%%%%%%%%%%%%%%%%%%%%%%%%%%%%%%%
\section{Uniform Distribution}
\label{app:A}
\label{app:uniform}
%%%%%%%%%%%%%%%%%%%%%%%%%%%%%%%%%%%%%%%%%%%%%%%%%%%%%%%%%%%%%%%%%%%%%%%%%%

In this appendix, we sketch the case of the uniform distribution on $\Sspace^n$.
The sole difference to the Poisson point process is that the
number of points is prescribed rather than a random variable.
Setting this number to $N = \density \sigma_{n+1}$,
it makes sense that in the limit,
when $N$ and $\density$ go to infinity,
the expected numbers of intervals of the radius function are the same
under both probabilistic models.
This is indeed what we establish now more formally.
By linearity of expectation, the number of intervals of type $(\ell, k)$
and geodesic radius at most $\GRad_0$ is
\begin{align}
  \Expected{\ccon{\ell}{k}{n}, \GRad_0}
    &=  \binom{N}{k+1} \Expected{\Pempty{\xxx}
                                 \cdot \One{k-\ell} (\xxx) 
                                 \cdot \One{\GRad_0} (\xxx)},
  \label{eqn:SMU_0}
\end{align}
in which $\xxx = (x_0, x_1, \ldots, x_k)$ is a sequence of $k+1$ points
on $\Sspace^n$,
$\GRad$ is the geodesic radius of the smallest circumscribed cap of $\xxx$,
and $\Pempty{\xxx} = ( 1 - \Area{\GRad}/\sigma_{n+1} )^{N-k+1}$
is the probability that this cap is empty.
The analogue of \eqref{eqn:SM} is therefore
\begin{align}
  \Expected{\ccon{\ell}{k}{n}, \GRad_0}
    &=  \binom{N}{k+1} \tfrac{1}{\sigma_{n+1}^{k+1}}
     \int\displaylimits_{\xxx \in (\Sspace^n)^{k+1}}
       \Pempty{\xxx} \cdot \One{k-\ell} (\xxx)
                     \cdot \One{\GRad_0} (\xxx) \diff \xxx .
  \label{eqn:SMU}
\end{align}
We apply the rotation-invariant Blaschke--Petkantschin formula
\eqref{eqn:rotationally-symmetric},
again with narrow bump functions as in \eqref{eqn:tryit}.
This gives
\begin{align}
  \!\! \Expected{\ccon{\ell}{k}{n}, \GRad_0}
    &=  \tfrac{N!}{(N-k-1)! \sigma_{n+1}^k}
        \tfrac{\sigma_n^k}{2 \Gama{k} n^{k-1}}
        \cdot \Ccon{\ell}{k}{n}
        \!\! \int\displaylimits_{t = 0}^{\sin^2 \GRad_0}
           \!\!\! t^{\frac{kn-2}{2}} (1\!-\!t)^{\frac{n-k-1}{2}}
           \left( \! 1 \!-\! \tfrac{\Area{\GRad}}{\sigma_{n+1}}
                                         \! \right)^{N-k+1} \!\! \diff t , 
  \label{eqn:f10_u}		
\end{align}
in which $\GRad = \GRad(t) = \arcsin \sqrt{t}$;
compare with \eqref{eqn:f10last}.
To prepare the next step, we note that
\begin{align}
  ( 1 - \tfrac{\Area{\GRad (t)}}{\sigma_{n+1}} )^{N-k+1}
    &\approx  e^{-\frac{N}{\sigma_{n+1}} \Area{\GRad(t)}}
\end{align}
as $t \to 0$.
From here on, we retrace the steps we took from \eqref{eqn:intuit1}
to \eqref{eqn:intuit2}.
In particular, we change variables with
$\tau = \tfrac{N}{\sigma_{n+1}} \nu_n t^{n/2}$,
and we substitute $\NRad_0 {\density^{-1/n}}$ for $\GRad_0$.
Observing $\tfrac{N!}{(N-k-1)!} \approx N^{k+1}$,
we simplify the expression and get
\begin{align}
  \Expected{\ccon{\ell}{k}{n}, \NRad_0}
    &=  N \cdot \tfrac{\iGama{v}{k}}{\Gama{k}} \cdot \Ccon{\ell}{k}{n}
        + \littleoh{N}
\end{align}
for the expected number of intervals of the radius function of the
Delaunay mosaic for $N$ points chosen uniformly at random on $\Sspace^n$,
in which $v = \NRad_0^n \nu_{n}$.
Comparing with the asymptotic result \eqref{eqn:asymp_res} in
Theorem \ref{thm:MainResult},
we see the same constants as for the Poisson point process.
However, the variance distinguishes the two cases,
being smaller for the uniform distribution than for the Poisson point
process; see \cite{Ste14}.

\end{document}